\documentstyle{article}
\font\germ=eufm10
\def\ssl{\hbox{\germ sl}}
\def\slh{\widehat{\ssl_2}}
\makeatletter
\def\aaa{@}
\title{\huge Crystallized structure  for level 0 part
of modified quantum affine algebra}
\author{Toshiki Nakashima
\thanks{supported by
the Ministry of Education, Science and Culture of Japan
as a overseas research scholar.
Current address : Department of Mathematics,
Northeastern University, Boston, MA 02115, USA.
(e-mail : toshiki\aaa neu.edu)}
\\
Department of Mathematical Science, Faculty of Engineering Science,\\
Osaka University, Toyonaka, Osaka 560, Japan\\
toshiki\aaa sigmath.es.osaka-u.ac.jp}
\date{}
\makeatother
  
\begin{document}
  \maketitle

   \renewcommand{\labelenumi}{$(\roman{enumi})$}
  \font\germ=eufm10

  \def\aff{{\rm Aff}}
  \def\al{\alpha}
  \def\beq{\begin{equation}}
  \def\beqn{\begin{eqnarray}}
  \def\beqnn{\begin{eqnarray*}}
  \def\bigsl{{\hbox{\fontD \char'54}}}
  \def\binf{B_{\infty}}
  \def\bminf{B_{-\infty}}
  \def\bl{\bullet}
  \def\bmax{B^{\scriptstyle{\rm max}}}
  \def\catob{{\cal O}(B)}
  \def\cd{\cdots}
  \def\del{\delta}
  \def\Del{\Delta}
  \def\ei{e_i}
  \def\eit{\tilde{e}_i}
  \def\ep{\epsilon}
  \def\eeq{\end{equation}}
  \def\eeqn{\end{eqnarray}}
  \def\eeqnn{\end{eqnarray*}}
  \def\FF{\hbox{\bf F}}
  \def\fit{\tilde{f}_i}
  \def\ge{\hbox{\germ g}}
  \def\gl{\hbox{\germ gl}}
  \def\gc{\bigcirc}
  \def\hom{{\hbox{Hom}}}
  \def\io{\iota}
  \def\llra{\relbar\joinrel\relbar\joinrel\relbar\joinrel\rightarrow}
  \def\lan{\langle}
  \def\lar{\longrightarrow}
  \def\lm{\lambda}
  \def\Lm{\Lambda}
  \def\mapright#1{\smash{\mathop{\longrightarrow}\limits^{#1}}}
  \def\mapleftright#1{\smash{\mathop{\longleftrightarrow}\limits^{#1}}}
  \def\map#1{\smash{\mathop{\longmapsto}\limits^{#1}}}
  \def\mpath{{\cal P}_m}
  \def\mlpath{{\cal P}_{m,l}}
  \def\mpathtc{{\cal P}_m(n;\vec{t};\vec{c})}
  \def\mlpathtc{{\cal P}_{m,l}(n;\vec{t};\vec{c})}
  \def\mllpathtc{{\cal P}_{m,l'}(n;\vec{t};\vec{c})}
  \def\nd{\noindent}
  \def\ot{\otimes}
  \def\op{\oplus}
  \def\ovl{\overline}
  \def\qq{\qquad}
  \def\q{\quad}
  \def\qed{\hfill\framebox[2mm]{}}
  \def\QQ{\hbox{\bf Q}}
  \def\qi{q_i}
  \def\qii{q_i^{-1}}
  \def\ran{\rangle}
  \def\SFF{\scriptstyle\hbox{\bf F}}
  \def\ssl{\hbox{\germ sl}}
  \def\slh{\widehat{\ssl_2}}
  \def\syl{\scriptstyle}
  \def\ti{t_i}
  \def\tii{t_i^{-1}}
  \def\til{\tilde}
  \def\tt{{\hbox{\germ{t}}}}
  \def\ttt{\hbox{\germ t}}
  \def\uq{U_q(\ge)}
  \def\uqm{U^-_q(\ge)}
  \def\uqmq{{U^-_q(\ge)}_{\bf Q}}
  \def\uqq{U^{\bf Q}_q(\ge)}
  \def\util{\widetilde\uq}
  \def\vmax{V^{\scriptstyle{\rm max}}}
  \def\vep{\varepsilon}
  \def\vp{\varphi}
  \def\vpi{\varphi^{-1}}
  \def\wtil{\widetilde}
  \def\ZZ{\hbox{\bf Z}}
\vspace{-26pt}
\begin{abstract}
\nd Crystal base of the level 0 part of
the modified quantum affine algebra
$\widetilde U_q(\widehat{sl_2})_0$  is given by path.
Weyl group actions, extremal vectors and
crystal structure of all irreducible components
 are described explicitly.
\end{abstract}

\section{Introduction}
\setcounter{equation}{0}
\renewcommand{\theequation}{\thesection.\arabic{equation}}

The modified quantum algebra, which is denoted $\util$,
was introduced in \cite{BLM} for $GL_n$-case
and in \cite{K2} for general case.
In \cite{L2},
G.Lusztig showed the existence of canonical (crystal) base
of modified quantum algebras for general Lie algebra.

In \cite{K4}, M.Kashiwara described detailed crystal structure of
the modified quantum algebras, in particular, he gave the
Peter-Weyl type decomposition theorem for the cases that
$\ge$ is finite type and
affine type with non-zero level(=central charge) parts. But, in
\cite{K4}, \cite{K5}, it is mentioned that
 the structure of level 0 part for affine type
is still unclear.
By the definition of the modified quantum algebra (\ref{ualm}) and
(\ref{eqn:def-util}),
we know that originally $\util$ is neither a highest nor
a lowest weight module. Nevertheless,
 if $\ge$ is affine, we can apply the powerful tool :
theory of integrable highest (resp. lowest) weight modules
to the positive (resp. negative)
level part
$\util_{+}:=\oplus_{\lan c,\lm\ran>0}\uq a_{\lm}$
(resp. $\util_{-}:=\oplus_{\lan c,\lm\ran<0}\uq a_{\lm}$)
by virtue of Weyl group actions on crystal bases,
where $c$ is a canonical central element of $\ge$. But,
in the level 0 case, there is no such a tool. However,
even in level 0 case,
it is still a good way to consider Weyl group actions
on crystal bases. Classification of 'extremal vectors'
(Definition \ref{extremal-vector}) is a crucial point in this paper.
By applying this classification to "path" realization, we can clarify
crystallized structure  of the level 0 part of the modified
quantum algbra for $\ge =\slh$ case and
give an explicit description of
its every connected component as a crystal graph.
The Peter-Weyl type theorem for this case
will be given in the forthcoming paper.

The path realization for the level 0 part of
the modified quantum algebra for $\ge =\slh$ case
has an another feature, which is a physical one.
A set of "path" is like the following thing:
\beq
\left\{(\cdots, i_k,i_{k+1},\cdots,)\,;\,
\hspace{-10pt}
\begin{array}{ll}
&i_k\in\ZZ,\,i_k=0\,(k\ll0),\\
&i_k=-i_{k+1}\,(k\gg0).
\end{array}
\right\}.
\label{(1)}
\eeq
Meanwhile, there is  so called 'XXZ type chain model',
which is a kind of physical model on the following space:
$$
{\hbox{\germ F}}=(\cdots\ot {\hbox{\bf C}}^{l+1}\ot
{\hbox{\bf C}}^{l+1}\ot{\hbox{\bf C}}^{l+1}\ot\cdots)^*,
$$
where ${\hbox{\bf C}}^{l+1}$ has a basis $\{(i)\}_{i=0,\cdots,l}$
and the notation $(\cdots)^*$ implies the condition that
${\hbox{\germ F}}$ is spanned by vectors
$\cdots\ot (i_k)\ot(i_{k+1})\ot\cdots$ with
$i_k+i_{k+1}=l$ for $|k|\gg0$.
We can see that this condition is similar to
the condition in (\ref{(1)}).
It is known that the space ${\hbox{\germ F}}$
has a $U_q(\slh)$-module structure.
In fact, in \cite{DFJMN} and \cite{IIJMNT},
this space is realized as
$$
{\hbox{\germ F}}
=\bigoplus_{\lan c,\zeta\ran=\lan c,\mu\ran=l}V(\lm)\ot V(-\mu),
$$
where $V(\zeta)$ (resp. $V(-\mu)$) is
an integrable highest (resp. lowest) weight module.
By considering the map $\pi_{\zeta,\mu}$ in (\ref{pi}),
we  can deduce that $\uq a_{\lm}$ is a kind of limit of
${\hbox{\germ F}}$ and
Theorem \ref{lusz} guarantees that
such a deduction is valid in the crystallized
space.
 For such a limit, in \cite{T} we gave
some related algebra structure and its representation theory.
But in this paper, we do not touch this subject any more.

Let us see  the organization of this paper.
In Sec.2,
we shall give the definitions of quantized enveloping algebra and crystal.
In Sec.3,
we introduce modified quantum algebra, its crystal base and Wely group
actions. From Sec.4 to Sec.7,
we consider the specific case $\ge=\slh$.
In Sec.4,
we study affinization of classical crystal and give a classification
of extremal vectors in $B^{\ot n}$,
where $B=\{\pm\}$ is the two-dimensional crystal, called 'spin'.
In Sec.5,
we shall give "path realization" of $\uq a_{\lm}$
with level of $\lm=0$ and introduce notions of
'domain' and 'wall', which play
a crucial role in this paper.
We also describe the actions of $\eit$ and $\fit$ on a path.
In Sec.6,
we give the path-spin correspondence,
which is a morphism of classical crystal between paths and spins.
In Sec.7,
first of all, we shall introduce some parametrizations
which are necessary to describe
connected components in $B(\util)$. Then
we shall give explicit crystallized structure of $\util$ by
classifying all extremal vectors in $\util$.

The author would like to acknowledge professor Masaki Kashiwara
for his helpful advises.
This work was partly done during the stay of the author
at Northeastern University.
He is grateful to professor Andrei Zelevinsky
for his kind hospitality.

%%%%%% Section 2 %%%%%%
\section{Notations and Preliminaries}
\setcounter{equation}{0}
\renewcommand{\theequation}{\thesection.\arabic{equation}}

\subsection{Definition of $U_q(\ge)$}

We shall recall the definition of
the quantized universal enveloping algebra.
First, let $\ge$ be
a  symmetrizable Kac-Moody algebra over {\bf Q}
with a Cartan subalgebra
$\ttt$, $\{\al_i\in\ttt^*\}_{i\in I}$
 the set of simple roots and
$\{h_i\in\ttt\}_{i\in I}$  the set of coroots,
where $I$ is a finite index set. We define an inner product on
$\ttt^*$ such that $(\al_i,\al_i)\in{\bf Z}_{\geq 0}$ and
$\lan h_i,\lm\ran=2(\al_i,\lm)/(\al_i,\al_i)$
for $\lm\in\ttt^*$.
Set $Q=\oplus_i\ZZ\al_i$,
$Q_+=\oplus_i\ZZ_{\geq0}\al_i$ and
$Q_-=-Q_+$. We call $Q$ a root lattice.
Let  $P$  a lattice of $\ttt^*$ {\it i.e.} a free
{\bf Z}-submodule of $\ttt^*$ such that
$\ttt^*\cong {\hbox{\bf Q}}\ot_{\ZZ}P$,
and $P^*=\{h\in \ttt|\lan h,P\ran\subset\ZZ\}$.
We set $P_+=\{\lm\in P|\lan \lm,h_i\ran\geq 0
{\hbox{ for any }}i\in I\}$.
An element of $P$(resp.$P_+$) is called a integral weight
(resp. dominant integral weight).

The quantized enveloping algebra $\uq$ is an associative
$\QQ(q)$-algebra generated by $\ei$, $f_i(i\in I)$
and $q^h(h\in P^*)$
satisfying the following relations:
\begin{eqnarray}
&&q^0=1, \q{\hbox{\rm and }}\q q^hq^{h'}=q^{h+h'},\\
&&q^he_iq^{-h}=q^{\lan h,\al_i\ran}e_i,\qq
q^hf_iq^{-h}=q^{-\lan h,\al_i\ran}f_i,\\
&&[e_i,f_j]=\del_{i,j}(t_i-t^{-1}_i)/(q_i-q^{-1}_i),\\
&&\sum_{k=1}^{1-{\lan h_i,\al_j\ran}}
(-1)^kx_i^{(k)}x_jx_i^{(1-{\lan h_i,\al_j\ran}-k)}=0,\qq(i\ne j)
\end{eqnarray}
where $x=e,f$ and we set $q_i=q^{(\al_i,\al_i)/2}$,
$t_i=q_i^{h_i}$, $[n]_i=(q^n_i-q^{-n}_i)/(q_i-q_i^{-1})$,
$[n]_i!=\prod_{k=1}^n[k]_i$,
$e_i^{(n)}=e_i^n/[n]_i!$ and $f_i^{(n)}=f_i^n/[n]_i!$.

It is well-known that $\uq$ has a Hopf algebra structure with a
comultiplication $\Del$ given by
$$
\Del(q^h)=q^h\ot q^h,\q \Del(\ei)=\ei\ot\tii+1\ot\ei,\q
\Del(f_i)=f_i\ot 1+\ti\ot f_i,
$$
for any $i\in I$ and $h\in P^*$.
We do not describe an antipode and a counit.
By this comultiplication, a tensor product of $\uq$-modules has
a $\uq$-module structure.

\subsection{Crystals}
\label{crystals}

Let us recall the definition of crystals \cite{K3,K4}.
The notion of a crystal
is motivated by abstracting the some combinatorial properties
of crystal bases.
We do not write down the definition of crystal base here
(See {\it e.g.} \cite{K1,K5,KN}).

\newtheorem{df}{Definition}[section]
\begin{df}
A {\it crystal} $B$ is a set with the following data:
\begin{eqnarray}
&&{\hbox{a map}}\q wt:B\lar P,\\
&&\vep_i:B\lar\ZZ\sqcup\{-\infty\},\q
  \vp_i:B\lar\ZZ\sqcup\{-\infty\},\q{\hbox{for}}\q i\in I,\\
&&\eit:B\lar B\sqcup\{0\},
\q\fit:B\lar B\sqcup\{0\}\q{\hbox{for}}\q i\in I.
\end{eqnarray}
Here 0 is an ideal element which is not included in $B$.
They are subject to the following axioms: For $b$,$b_1$,$b_2\in B$,
\begin{eqnarray}
&&\vp_i(b)=\vep_i(b)+\lan h_i,wt(b)\ran,\\
&&wt(\eit b)=wt(b)+\al_i{\hbox{ if }}\eit b\in B,\\
&&wt(\fit b)=wt(b)-\al_i{\hbox{ if }}\fit b\in B,\\
&&\eit b_2=b_1 {\hbox{ if and only if }} \fit b_1=b_2,
\label{eeff}\\
&&{\hbox{if }}\vep_i(b)=-\infty,
  {\hbox{ then }}\eit b=\fit b=0.
\end{eqnarray}
\end{df}
{}From the axiom (\ref{eeff}),
we can consider the graph strucure on a crystal $B$.
\begin{df}
\label{c-gra}
The crystal graph of crystal $B$ is
an oriented and colored graph given by
the rule : $b_1\mapright{i} b_2$ if and only if $b_2=\fit b_1$
$(b_1,b_2\in B)$.
\end{df}
\begin{df}
\label{df:mor}
\begin{enumerate}
\item
If $B$ has the weight decomposition $B=\bigsqcup_{\lm\in P}B_{\lm}$
where $B_{\lm}=\{b\in B|wt(b)=\lm\}$ for $\lm \in P$, we call
$B$ a $P$-weighted crystal.
\item
Let $B_1$ and $B_2$ be crystals.
A morphism of crystals $\psi:B_1\lar B_2$
is a map $\psi:B_1\sqcup\{0\}\lar B_2\sqcup\{0\}$
satisfying the following axioms:
\begin{eqnarray}
&&\hspace{-30pt}\psi(0)=0,
\label{psi(0)=0}\\
&&\hspace{-30pt}wt(b)=wt(\psi(b)),\q \vep_i(b)=\vep_i(\psi(b)),\q
\vp_i(b)=\vp_i(\psi(b))
\label{wt}\\
&&{\hbox{if }}b\in B_1{\hbox{ and }}\psi(b)\in B_2,\nonumber\\
&&\hspace{-30pt}\psi(\eit b)
=\eit\psi(b){\hbox{ if }}b\in B_1{\hbox{ satisfies }}
 \psi(b)\neq0{\hbox{ and }}\psi(\eit b)\neq0,\\
&&\hspace{-30pt}\psi(\fit b)
=\eit\psi(b){\hbox{ if }}b\in B_1{\hbox{ satisfies }}
 \psi(b)\neq0{\hbox{ and }}\psi(\fit b)\neq0.
\end{eqnarray}
\item
A morphism of crystals $\psi:B_1\lar B_2$ is called {\it strict} if the
associated map from $B_1\sqcup\{0\}\lar B_2\sqcup\{0\}$ commutes with
all $\eit$ and $\fit$.
If $\psi$ is injective, surjective and strict,
$\psi$ is called an {\it isomorphism}.
\item
A crystal $B$ is a {\it normal}, if for any subset $J$ of $I$ such that
$((\al_i,\al_j))_{i,j\in J}$ is a positive symmetric matrix,
$B$ is isomorphic
to a crystal base of an integrable $U_q(\ge_{J})$-module, where
$U_q(\ge_{J})$ is the quantum algebra generated by $e_j$,
$f_j$ $(j\in J)$ and $q^h$ $(h\in P^*)$.
\end{enumerate}
\end{df}

For crystals $B_1$ and $B_2$, we shall define their tensor product
$B_1\ot B_2$ as follows:
\begin{eqnarray}
&&B_1\ot B_2=\{b_1\ot b_2| b_1\in B_1 ,\, b_2\in B_2\},\\
&&wt(b_1\ot b_2)=wt(b_1)+wt(b_2),\\
&&\vep_i(b_1\ot b_2)={\hbox{max}}(\vep_i(b_1),
  \vep_i(b_2)-\lan h_i,wt(b_1)\ran),
\label{tensor-vep}\\
&&\vp_i(b_1\ot b_2)={\hbox{max}}(\vp_i(b_2),
  \vp_i(b_1)+\lan h_i,wt(b_2)\ran),
\label{tensor-vp}\\
&&\eit(b_1\ot b_2)=
\left\{
\begin{array}{ll}
\eit b_1\ot b_2 & {\mbox{ if }}\vp_i(b_1)\geq \vep_i(b_2)\\
b_1\ot\eit b_2  & {\mbox{ if }}\vp_i(b_1)< \vep_i(b_2),
\end{array}
\right.
\label{tensor-e}
\\
&&\fit(b_1\ot b_2)=
\left\{
\begin{array}{ll}
\fit b_1\ot b_2 & {\mbox{ if }}\vp_i(b_1)>\vep_i(b_2)\\
b_1\ot\fit b_2  & {\mbox{ if }}\vp_i(b_1)\leq \vep_i(b_2),
\label{tensor-f}
\end{array}
\right.
\end{eqnarray}

Here we understand that $0\ot b=b\ot 0=0$.
Let ${\cal C}(I,P)$ be the category
of crystals determined by the data $I$ and $P$.
Then $\ot$ is a functor
from ${\cal C}(I,P)\times{\cal C}(I,P)$
to ${\cal C}(I,P)$ and satisfies the
associative law:
$(B_1\ot B_2)\ot B_3\cong B_1\ot(B_2\ot B_3)$ by
$(b_1\ot b_2)\ot b_3\leftrightarrow b_1\ot (b_2\ot b_3)$.
 Therefore, the category of crystals is endowed
with the structure of tensor category.

\newtheorem{ex}[df]{Example}
\begin{ex}\label{Example:crystal}
We give some examples of crystals.
\begin{enumerate}
\item  $C=\{c\}$ with
$$
wt(c)=0,\q \vep_i(c)=\vp_i(c)=0,\q \eit c=\fit c=0.
$$
This is isomorphic to the crystal of the trivial $\uq$-module.
\item For $\lm\in P$, we set $T_{\lm}=\{t_{\lm}\}$ with
$$
wt(t_{\lm})=\lm,\q \vep_i(t_{\lm})=\vp_i(t_{\lm})=-\infty,\q
\eit(t_{\lm})=\fit(t_{\lm})=0.
$$
We can see that $T_{\lm}\ot T_{\mu}\cong T_{\lm+\mu}$ and
$B\ot T_{0}\cong T_{0}\ot B\cong B$ for any crystal $B$.
\item
For $\lm\in P_+$, let $(L(\lm),B(\lm))$ be the crystal base of a
$\uq$-integrable highest weight module $V(\lm)$  $($\cite{K3}$)$.
$B(\lm)$ is the crystal associated with $V(\lm)$.
Let $B(-\lm)$ be the crystal associated with a integrable lowest
weight module $V(-\lm)$. Then $C$ is isomorphic to $B(0)$.
\item Let $(L(\infty),B(\infty))$ be a crystal base of $\uqm$
$($\cite{K3}$)$.
$B(\infty)$ is a crystal associated with $\uqm$.
We shall also denote $B(-\infty)$ for a crystal associated
with $U^+_q(\ge)$.
\end{enumerate}
\end{ex}

%%%%%%%%% This section is Section 3 %%%%%%%%
\section{Crystals of modified quantum algebra}
\setcounter{equation}{0}
\renewcommand{\theequation}{\thesection.\arabic{equation}}
This section is devoted to
review \cite{K4},\cite{L2} (See also \cite{L1}).
\subsection{Modified quantum algebra and Crystal base}

For an integral weight $\lm\in P$,
let $\uq a_{\lm}$ be the left $\uq$-module
given by
\beq
\uq a_{\lm}:=\uq/\sum_{h\in P^*}\uq(q^h-q^{\lan h,\lm\ran}),
\label{ualm}
\eeq
where $a_{\lm}$ is the image of the unit by the canonical projection.
We set
\begin{equation}
\widetilde\uq=\bigoplus_{\lm\in P}\uq a_{\lm},
\label{eqn:def-util}
\end{equation}
which is called {\it modified quantum algebra}.

We shall see a crystal base of $\widetilde\uq$.
Taking $\lm\in P$ and
 choosing $\zeta,\mu\in P_+$ such that $\lm=\zeta-\mu$,
we get the following $\uq$-linear surjective homomorphism:
\beqn
\pi_{\zeta,\mu}  : \uq a_{\lm}&\longrightarrow &V(\zeta)\ot V(-\mu),
\label{pi}\\
               a_{\lm}&\mapsto &u_{\zeta}\ot u_{-\mu}.\nonumber
\eeqn
where $V(\zeta)$ and $V(-\mu)$ are
as in Example \ref{Example:crystal} (iii) and
$u_{\zeta}$ and $u_{-\mu}$ are their highest weight vector
and lowest weight vector respectively.

\newtheorem{thm3}{Theorem}[section]
\begin{thm3}[{\it cf} \cite{L2}]
\label{lusz}
For any $\lm\in P$,
there exists a unique pair $(L(\uq a_{\lm}),B(\uq a_{\lm}))$
which satisfies the following properties.
\begin{enumerate}
\item
We set $A:=\{f(q)\in\QQ(q)|{\hbox{$f$ has no pole at $q=0$}}\}$.
$L(\uq a_{\lm})$ is a free $A$-module
such that $\uq a_{\lm}\cong \QQ(q)\ot_A L(\uq a_{\lm})$ and
$B(\uq a_{\lm})$ is a $\QQ$-basis of
the $\QQ$-vector space $L(\uq a_{\lm})/qL(\uq a_{\lm})$.
\item
For any $\zeta,\mu\in P_+$ with $\lm=\zeta-\mu$, we have
$$
\pi_{\zeta,\mu} (L(\uq a_{\lm}))\subset L(\zeta)\ot_A L(-\mu),
$$
and the induced map $\bar\pi_{\zeta,\mu}$:
$$
\bar\pi_{\zeta,\mu}\,:\,L(\uq a_{\lm})/qL(\uq a_{\lm})\longrightarrow
(L(\zeta)/qL(\zeta))\ot(L(-\mu)/qL(-\mu)),
$$
satisfies $\bar\pi_{\zeta,\mu}
(B(\uq a_{\lm}))\subset B(\zeta)\ot B(-\mu)\,\sqcup \,\{0\}$.
\item
There is a structure of crystal
on $B(\uq a_{\lm})$ such that $\bar\pi_{\zeta,\mu}$ gives
a strict morphism of crystals
for any $\zeta,\mu\in P_+$ with $\lm=\zeta-\mu$.
\end{enumerate}
\end{thm3}

\vskip7pt
\nd
Set
$$
(L(\util),B(\util)):=\bigoplus_{\lm\in P}(L(\uq a_{\lm}),B(\uq a_{\lm})).
$$
{\sl Remark.}
$B(\uq a_{\lm})$ is a normal crystal and then $B(\util)$ is a normal crystal.

Let $B(\infty)$, $B(-\infty)$ and $T_{\lm}$ $(\lm\in P)$
be the crystals given in Example \ref{Example:crystal}.
The following theorem plays a significant role in this paper
(See \cite[Sec.3]{K4}).

\begin{thm3}
\label{U=BTB}
$B(\uq a_{\lm})\cong B(\infty)\ot T_{\lm}\ot B(-\infty)$ as a crystal.
\end{thm3}

\newtheorem{cor3}[thm3]{Corollary}
\begin{cor3}
$B(\wtil\uq)\cong \oplus_{\lm\in P}B(\infty)\ot T_{\lm}\ot B(-\infty)$
as a crystal.
\end{cor3}

\subsection{Weyl group action and Extremal vectors}

This subsection is devoted to review \cite[Sec.7.8.9]{K4}.
Let $B$ be a normal crystal (See \ref{crystals}).
 Let us define the Weyl group action
on the underlying set $B$. For $i\in I$ and $b\in B$, we set
\begin{equation}
S_ib=
\left\{
\begin{array}{ll}
\fit^{\lan h_i,wt(b)\ran}b &{\rm if }\,\, \lan h_i,wt(b)\ran\geq 0\\
\eit^{-\lan h_i,wt(b)\ran}b &{\rm if }\,\, \lan h_i,wt(b)\ran<0.
\end{array}
\right.
\end{equation}

We can easily obtain the following formula:
$$
 S_i^2={\hbox{id}},\qq
 S_i\eit=\fit S_i,\qq
 wt(S_i b)=s_i(wt(b)),
$$
where $s_i(\lm)=\lm-\lan h_i,\lm\ran\al_i$ is the simple reflection.

Let $\ge$ be a rank 2 finite dimensional Lie algebra,
and $W$ be the Weyl group associated with $\ge$.
We set $w_0=s_{i_1}\cdots s_{i_k}$ a reduced expression of
the longest element of $W$. Here we get the following (\cite[Sec.7]{K4}):

\newtheorem{pro3}[thm3]{Proposition}
\begin{pro3}
Let $B$ be a normal crystal.
For any $b\in B$, $S_{i_1}\cdots S_{i_k}b$ does not depend on the choice of
reduced expression.
\end{pro3}

\begin{cor3}
$\{S_i\}_{i=1,2}$ satisfies the braid relation.
\end{cor3}

Thus for general $\ge$,
we know that $\{S_i\}_{i\in I}$ defines the Weyl group
action on a normal crystal $B$.
\newtheorem{df3}[thm3]{Definition}
\begin{df3}
\label{extremal-vector}
\begin{enumerate}
\item  Let $B$ be a normal crystal.
An element $b\in B$ is called $i$-extremal,
       if $\eit b=0$ or $\fit b=0$.
\item  An element $b\in B$ is called extremal if for any $l\geq0$,
$S_{i_1}\cdots S_{i_l}b$ is $i$-extremal for any $i$,
$i_1\cdots i_l\in I$.
\end{enumerate}
\end{df3}

The following theorems play a significant role in Sec.7.
\begin{thm3}
\label{conn-ext}
Any connected component of $B(\util)$ contains an extremal vector.
\end{thm3}

%%%%%%%%%%%%%%%%%%%% This section is Section 4 %%%%%%%%%%%%
\section{Affine crystals}
\setcounter{equation}{0}
\renewcommand{\theequation}{\thesection.\arabic{equation}}
In the rest of this paper, we fix $\ge=\widehat{\ssl_2}$.
Let us denote $U$ for $U_q(\slh)$.
\subsection{Notations}
We follow the notations in  \cite{IIJMNT}, \cite{KMN}.
Let  $\Lm_0$ and $\Lm_1$ be fundamental weights and
$\del$ a generator of null roots. Then
we can write a weight lattice
$P=\ZZ\Lm_0\op \ZZ\Lm_1\op\ZZ\del$, its
dual lattice $P^*=\ZZ h_0\op\ZZ h_1\op \ZZ d$ and the Cartan subalgebra
$\tt=\QQ h_0\op\QQ h_1\op \QQ d$, where
$\lan h_i,\Lm_j\ran=\del_{ij}$, $\lan h_i,\del\ran=\lan d,\Lm_i\ran=0$,
$\lan d,\del\ran=1$ $\al_1=2\Lm_1-2\Lm_0$ and $\al_0=\del-\al_1$.
Set $P_{cl}=P/\ZZ\del$ and
let $cl:P\longrightarrow P_{cl}$ be a canonical projection.
Then we set $(P_{cl})^*=\ZZ h_0\op\ZZ h_1$.
Now, we fix a map $af:P_{cl}\longrightarrow P$
such that $af\circ cl(\Lm_i)=\Lm_i$.
Here note that $af\circ cl(\al_1)=\al_1$ and
$al\circ cl(\al_0)=\al_0-\del$.
In this specific case,
we call an element of $P$ an {\it affine weight}
and an element of $P_{cl}$
a {\it classical weight}.  In the rest of this paper,
if there is no confusion, we shall denote $\Lm_i$ for
a classical weight $cl(\Lm_i)$.

\subsection{Affinization of classical crystals}
\label{4.2}
$U$ is the quantized enveloping algebra associated with $P$.
Let $U'$ be its subalgebra generated by
$e_i$, $f_i$ and $q^h\q(h \in (P_{cl})^*)$.
The algebra $U'$ is also a
quantized enveloping algebra associated with $P_{cl}$.
Now, we call a $P$-weighted crystal an {\it affine crystal} and
a $P_{cl}$-weighted crystal a {\it classical crystal}.

\vskip6pt
\nd
{\sl Remark.}
\begin{enumerate}
\item
The simple roots of $U'$ : $\{cl(\al_i)\}$ are not linearly independent.
\item
A $U$-module has a $U'$-module structure but
in general, the opposite case is false.
\end{enumerate}

\newtheorem{def4}{Definition}[section]
\newtheorem{pr4}[def4]{Proposition}
\begin{def4}
\label{def4:aff}
Let $B$ be a classical crystal.
We define the affine crystal ${\rm Aff}(B)$
associated with $B$ as follows;
\begin{equation}
\aff(B):=\{z^n\ot b\,|\,b\in B,\,n\in \ZZ\},
\end{equation}
where $z$ is an indeterminate.
We call $\aff(B)$ an affinization of $B$.
The actions by $\eit$ and $\fit$, and the data are given as follows:
\begin{equation}
\begin{array}{rcl}
&&\eit(z^n\ot b)=z^{n+\del_{i,0}}\ot \eit(b),\q
\fit(z^n\ot b)=z^{n-\del_{i,0}}\ot \fit(b),\\
&&\vep_i(z^n\ot b)=\vep_i(b),\q \vp_i(z^n\ot b)=\vp_i(b),\q
wt(z^n\ot b)=n\del+af(wt(b)).
\end{array}
\label{delta}
\end{equation}
\end{def4}

Here, note that even if a classical crystal $B$ is  connected
as a crystal graph, its affinization $\aff(B)$ is not necessarily connected.

\newtheorem{ex4}[def4]{Example}
\begin{ex4}
\label{ex4:B}
Let $B=\{+,-\}$ be a 2-dimensional classical crystal given by
\beq
\begin{array}{rcl}
&&\til e_0(+)=\til f_1(+)=-,\q \til e_1(-)=\til f_0(-)=+,
\\
&&\til e_0(-)=\til f_1(-)=0,\q \til e_1(+)=\til f_0(+)=0,
\\
&&\vp_1(+)=\vp_0(-)=\vep_0(+)=\vep_1(-)=1,\\
&&\vp_1(-)=\vp_0(+)=\vep_0(-)=\vep_1(+)=0,
\\
&& wt(\pm)=\pm(\Lm_1-\Lm_0).
\end{array}
\label{act}
\eeq
It is easy to see that $B^{\ot 2}$ is connected. But its affinization
\begin{equation}
\aff(B^{\ot 2})\cong
\{z^n\ot \epsilon_1\ot \epsilon_2|\epsilon_i=\pm,\q n\in \ZZ\}
\end{equation}
is not connected.
In fact, this is divided into the following two components;
\begin{eqnarray}
&&\aff(B^{\ot 2})_1:=\{z^{2n-1}\ot +\ot +,\,z^{2n-1}\ot -\ot +,\,
   z^{2n-1}\ot -\ot -,\,z^{2n}\ot +\ot -,|n\in \ZZ\}.\nonumber \\
&&\aff(B^{\ot 2})_0:=\{z^{2n}\ot +\ot +,\,z^{2n}\ot -\ot +,\,
   z^{2n}\ot -\ot -,\,z^{2n+1}\ot +\ot -,|n\in \ZZ\},\nonumber
\end{eqnarray}
\end{ex4}

\subsection{Extremal vectors in $B^{\ot n}$}
In this subsection, let $B$ be the 2-dimensional classical crystal
introduced in Example \ref{ex4:B}.

Now, we classify  all extremal vectors in $B^{\ot n}$.
\begin{pr4}
\label{ext}
$B^{\ot n}$ is connected as a crystal graph.
Let $E$ be the set of all extremal vectors in $B^{\ot n}$.
Then we have
$$
E=\{(+)^{\ot n}, (-)^{\ot n}\}.
$$
\end{pr4}

\noindent
{\sl Proof.}
By the fact that $B$ is a perfect crystal
\cite[Corollary 4.6.3.]{KMN},
we can easily obtain the connectedness of $B^{\ot n}$.

By \cite[Sec.2]{KN}, we know that
\beq
\begin{array}{rcl}
&&\til f_1((-)^{\ot k}\ot(+)^{\ot l})=(-)^{\ot k+1}\ot(+)^{\ot l-1},\\
&&\til e_1((-)^{\ot k}\ot(+)^{\ot l})=(-)^{\ot k-1}\ot(+)^{\ot l+1},\\
&&\til f_0((+)^{\ot k}\ot(-)^{\ot l})=(+)^{\ot k+1}\ot(-)^{\ot l-1},\\
&&\til e_0((+)^{\ot k}\ot(-)^{\ot l})=(+)^{\ot k-1}\ot(-)^{\ot l+1},
\end{array}
\label{4act}
\eeq
where we consider $(\pm)^{\ot m}=0$ if $m<0$.
Since $B^{\ot n}$ is a normal crystal, we get
\begin{equation}
S_1(+)^{\ot n}=S_0(+)^{\ot n}=(-)^{\ot n},\q
S_1(-)^{\ot n}=S_0(-)^{\ot n}=(+)^{\ot n}.
\label{eqn:SS}
\end{equation}
{}From (\ref{eqn:SS}) and the fact that
$\til e_1(+)^{\ot n}=\til e_0(-)^{\ot n}
=\til f_1(-)^{\ot n}=\til f_0(+)^{\ot n}=0
$, we know that $(+)^{\ot n}$ and $(-)^{\ot n}$ are extremal vectors.

Now we shall show that there is no other extremal vector without
these two vectors by using the induction on $n$.

For $n=1$, this is trivial. We assume that $u\ot +$ is an extremal vector
in $B^{\ot n+1}$.
Here note that
for any $b\in B^{\ot n}$, $\vp_i(b)$ and $\vep_i(b)$ are given by
\beqnn
&&\vp_i(b)={\rm max}\{n;\fit^nb\ne0\},\\
&&\vep_i(b)={\rm max}\{n;\eit^nb\ne0\},
\eeqnn
and then
\beq
\vp_i(b)\geq 0,\q \vep_i(b)\geq 0.
\label{vpvep>0}
\eeq
We have $\til f_1(u\ot +)\ne 0$
by $\til f_1(+)\ne 0$ and (\ref{tensor-f}).
Then, by the definition of extremal vectors, we have
$\til e_1(u\ot +)=0$ and then
\beq
\label{e1u=0}
\til e_1u=0,
\eeq
since  we have that $\vep_1(+)=0$, $\vp_1(u)\geq0$
and then $\til e_1(u\ot +)=\til e_1(u)\ot +$ by (\ref{tensor-e}).

\noindent
We shall show
\begin{equation}
\label{e0u+ne0}
\til e_0(u\ot +)\ne0.
\eeq
If $\vp_0(u)\geq1$,
$\vp_0(u)\geq 1=\vep_0(+)$ and
then $\til e_0(u\ot+)=\til e_0(u)\ot+\ne0$.
Otherwise, $\vp_0(u)<\vep_0(+)$ and then
$\til e_0(u\ot+)=u\ot\til e_0(+)=u\ot-\ne0$.
We get (\ref{e0u+ne0}).
By the definition of extremal vector, we have
\beq
\label{f0u+=0}
\til f_0(u\ot+)=0.
\eeq
By (\ref{e1u=0}), we have $\vep_1(u)=0$. Then we get
$\lan h_1, wt(u\ot+)\ran=\lan h_1, wt(u)\ran+1=\vp_1(u)+1>0$. Thus
\beqnn
S_1(u\ot +) & = & \til f_1^{\lan h_1, wt(u\ot+)\ran}(u\ot +)
 = \til f_1^{\vp_1(u)}u\ot \til f_1(+)  \\
& = & S_1(u)\ot-.
\eeqnn
This $S_1(u)\ot -$ is extremal and $\til e_1(u\ot+)\ne0$. Therefore,
by the similar argument as above, we get $\til e_0S_1(u)\ne0$ and
\beq
S_0(S_1(u)\ot-)=S_0S_1(u)\ot+.
\eeq
By arguing similarly, we get $\til e_1(S_0S_1u)=0$.
Repeating these argument, we obtain
\beq
\til e_1(S_0S_1\cdots S_1u) = 0,\,\,
\til e_0(S_1S_0\cdots S_1u) = 0.
\label{e1:01u=0}
\eeq
The set $\{(S_0S_1)^ku\}_{k\in \ZZ_{\geq0}}$ is
a subset of the finite set $B^{\ot n}$.
Then there exist $l,\,\,m\in \ZZ_{\geq0}$ such that $l>m$,
$$
(S_0S_1)^lu=(S_0S_1)^mu.
$$
Then we have $(S_0S_1)^{l-m}u=u$ ($l-m>0$)
and then for any $p\geq 0$ there exists $r\in\ZZ_{>0}$
such that $(l-m)r>p$.
Thus we have
\beq
(S_1S_0)^pu=(S_0S_1)^{(l-m)r-p}u,\,\,
S_0(S_1S_0)^{p}u=S_1(S_0S_1)^{(l-m)r-p-1}u.
\label{01=10}
\eeq
By (\ref{e1:01u=0}) and (\ref{01=10}),
we obtain
\beq
\til e_1(S_1S_0\cdots S_1S_0 )u=0,\,\,
\til e_0(S_0S_1\cdots S_1S_0)u=0.
\label{e1:10u=0}
\eeq

\noindent
We shall show
\beq
\til f_0u=0.
\label{f0u=0}
\eeq
Assuming $\til f_0u\ne 0$, we shall derive a contradiction.
The assumption implies $\vp_0(u)>0$. If $\vp_0(u)\geq2$,
$\til f_0(u\ot+)=\til f_0(u)\ot +\ne0$
since $\vep_0(+)=1$. This contradicts  (\ref{f0u+=0}).
Then we know that
\beq
\vp_0(u)=1.
\label{vp0u=1}
\eeq
Now we write $u=u_1\ot u_2\ot\cdots\ot u_n$ ($u_j=\pm$).
By using (\ref{e1u=0}) and (\ref{vp0u=1}), we get
\beqn
&& u_1=+,
\label{u1=+} \\
&&\vep_0(u)\geq \vp_0(u)=1,
\label{vep0>vp0}
\eeqn
because
$0\leq \lan h_1, wt(u)\ran=-\lan h_0,wt(u)\ran=\vep_0(u)-\vp_0(u)$.
Now, by applying Remark 2.1.2 in \cite{KN},
(\ref{u1=+}) and (\ref{vep0>vp0})
to $S_0\,u$. It is easily obtained that
$S_0u=\til e_0^{-\lan h_0,wt(u)\ran}u = \til e_0^{\vep_0(u)-1}u$
is in the following form:
\beq
S_0\,u=+\ot u',
\label{s0u=+u'}
\eeq
where $u'\in B^{\ot{n-1}}$, {\it i.e.}
the action of $S_0$ never touches $u_1$.
A vector in the form (\ref{s0u=+u'})
does not vanish by the action of $\til e_0$.
This contradicts  (\ref{e1:10u=0})
and then we get  $\vp_0(u)=0$.

\noindent
We have
\beqn
S_0(u\ot +) & = & \til e_0^{-\lan h_0,wt(u\ot +)\ran}(u\ot +)
            = \til e_0^{\vep_0(u)-\vp_0(u)+1}(u\ot+) \nonumber \\
            & = & \til e_0^{\vep_0(u)+1}(u\ot +)
            = \til e_0^{\vep_0(u)}\, u\ot \til e_0(+)\nonumber \\
            & = & S_0\,u\ot - \nonumber
\eeqn
The vector $S_0\,u\ot -$ does not vanish
by the action of $\til f_0$
since $\til f_0(-)\ne0$.
Since $S_0u\ot -$  is an extremal vector,
this vanishes by the action of $\til e_0$.
By the similar argument to obtain (\ref{e0u+ne0}),
we have $\til e_1(S_0\, u\ot -)\ne 0$.
Then
\beq
\til f_1(S_0\, u\ot -)=0.
\label{f1s0u-=0}
\eeq
By exchanging $+$ and $-$,
and arguing similarly to get (\ref{f0u=0}), we get
\beq
\til f_1(S_0\,u)=0.
\label{f1s0u=0}
\eeq
By repeating above argument, we obtain
\beq
\til f_0(S_1S_0\cdots S_0\, u)  = 0,\,\
\til f_1(S_0S_1\cdots S_0\, u)  = 0.
\label{f0:10u=0}
\eeq
Furthermore, by the similar argument to get (\ref{e1:10u=0}),
we get
\beqn
\til f_0(S_0S_1\cdots S_1\,u)  =  0,\,\,
\til f_1(S_1S_0\cdots S_1\,u)  =  0.
\label{f1:101u=0}
\eeqn
By (\ref{e1:01u=0}), (\ref{e1:10u=0}),
(\ref{f0:10u=0}) and (\ref{f1:101u=0}),
we know that the vector $u$ is an extremal vector
in $B^{\ot n}$. By the hypothesis of the induction and (\ref{u1=+}),
we get
$$
u=(+)^{\ot n} \,\,\,{\rm and \,\,then }\,\,\,u\ot +=(+)^{\ot{n+1}}.
$$
By assuming $u\ot -$ is an extremal vector in
$B^{\ot{n+1}}$ and discussing similarly,
we get $u=(-)^{\ot n}$ and then $u\ot -=(-)^{\ot{n+1}}$.
\hfill\framebox[2mm]{}

 %%%%%% section 5 %%%%%%%%%
\section{Path realization for $B(Ua_{\lm})$ with level of $\lm=0$}
\label{Sec5}
  \setcounter{equation}{0}
  \renewcommand{\theequation}{\thesection.\arabic{equation}}
  As in the previous section, we set $\ge=\slh$ in this
  section.

  \subsection{Crystal $B(\infty)$ and $B(-\infty)$}
  \label{5.1}
  Now, we define the following $\slh$-classical crystal:
  \newtheorem{def5}{Definition}[section]
  \begin{def5}
  We set
  $$
  B_{\infty}:=\{(n)|n\in \ZZ\}, \q(wt(n)=2n(\Lm_0-\Lm_1))
  $$
  and define the actions of $\eit$ and $\fit$ as follows;
  $$
  \til e_1(n)=(n-1),\,\til f_1(n)=(n+1),\,
  \til e_0(n)=(n+1),\,\til f_0(n)=(n-1).
  $$
  We also set
  $$
  \vep_1(n)=n,\,\vp_1(n)=-n,\,
  \vep_0(n)=-n,\,\vp_0(n)=n.
  $$
  \end{def5}
  By the above data, $B_{\infty}$ is equipped with a classical
  crystal structure.

  We introduce the following remarkable result (see
  \cite{KKM}).
  \newtheorem{pr5}[def5]{Proposition}
  \begin{pr5}
  Let $B_{\infty}$ be as above. We get the following
  isomorphism of classical crystal:
  \beqn
  \label{B=BB}
  B(\infty) &\mapright{\sim}& B(\infty)\ot B_{\infty}
  \q({\hbox{ resp. }}
  B(-\infty)\mapright{\sim} B_{\infty}\ot B(-\infty)),\\
  u_{\infty} &\mapsto& u_{\infty}\ot(0)\q
  ({\hbox{resp. }} u_{-\infty}\mapsto (0)\ot u_{-\infty}).\nonumber
  \eeqn
  \end{pr5}
  By applying this proposition repeatedly, we get for any
  $k>0$,
  \beqn
  \psi_k:B(\infty)& \mapright{\sim} & B(\infty)\ot B_{\infty}^{\ot k}\,
  ({\hbox{ resp. }}
  B(-\infty)\mapright{\sim} B_{\infty}^{\ot k}\ot B(-\infty)),
  \label{eqn:B-iso}\\
  u_{\infty}& \mapsto & u_{\infty}\ot (0)^{\ot k}\,
  ({\hbox{resp. }}
  u_{-\infty}\mapsto (0)^{\ot k}\ot u_{-\infty}).\nonumber
  \eeqn

\newtheorem{lem5}[def5]{Lemma}
  \begin{lem5}
  \label{lem5:path}
  For any $b\in B(\infty)$ (resp. $B(-\infty)$), there exists
  $k>0$ such that
  \begin{equation}
  \psi_k(b)\in u_{\infty}\ot B_{\infty}^{\ot k}
  \q({\hbox{ resp. }} B_{\infty}^{\ot k}\ot u_{-\infty}).
  \end{equation}
  \end{lem5}

\nd
 We set
  \begin{eqnarray}
  &&{\cal P}(\infty):=\{(\cdot\cdot,i_k,i_{k+1},\cdot\cdot,i_{-1})|
  i_k\in B_{\infty}{\hbox{ and }}if |k|\gg0, i_k=(0)\},\\
  &&{\cal P}(-\infty):=\{(i_0,\cdot\cdot,i_k,i_{k+1},\cdot\cdot)|
  i_k\in B_{\infty}{\hbox{ and }}if |k|\gg0, i_k=(0)\},
  \end{eqnarray}

\nd
Now, we consider formally $u_{\infty}=\cdots\ot  (0)\ot (0)=(\cdots,(0),(0))$
  (resp. $u_{-\infty}=(0)\ot (0)\ot \cdots=((0),(0),\cdots)$).
  Then by (\ref{eqn:B-iso}) and Lemma \ref{lem5:path},
  we get the following isomorphism between $B(\infty)$ (resp.
  $B(-\infty)$)
  and ${\cal P}(\infty)$ (resp. ${\cal P}(-\infty)$).
  \begin{pr5}
  \label{0-path}
  The crystal $B(\infty)$ $({\hbox{resp. }} B(-\infty))$ is isomorphic
  to
  ${\cal P}(\infty)$ $({\hbox{resp. }} {\cal P}(-\infty))$ given by
  $B(\infty)\ni b\mapsto p\in {\cal P}(\infty)$
  (resp. $B(-\infty)\ni b\mapsto p\in {\cal P}(-\infty)$)
  where $\psi_k(b)=u_{\infty}\ot i_k\ot\cdots\ot i_{-2}\ot
  i_{-1}$
  $({\hbox{resp. }} \psi_k(b)=i_0\ot i_1\ot\cdots\ot i_{k}\ot
  u_{-\infty})$
  for $|k|\gg0$.
  \end{pr5}

  \subsection{Path}

  \begin{lem5}
  We set $\lm=m(\Lm_0-\Lm_1)\in P_{cl}$ $(m\in \ZZ)$. Then the map
  \begin{eqnarray}
  \vp:& T_{\lm}\ot B_{\infty} & \mapright{\sim} B_{\infty}\ot
  T_{-\lm},
  \label{TB=BT}\\
      & t_{\lm}\ot (n) & \mapsto (m+n)\ot t_{-\lm},\nonumber
  \end{eqnarray}
  is an isomorphism between classical crystals.
  \end{lem5}

  \noindent
  {\sl Proof. }\,
  It is trivial that the map $\vp$ is bijective.
  Then we shall show that $\vp$ is a strict morphism of
  classical crystals
  (See Definition \ref{df:mor}).
  The weight of $t_{\lm}\ot (n)$ is $(m+2n)(\Lm_0-\Lm_1)$ and
  the one of
  $(m+n)\ot t_{-\lm}$ is also $(m+2n)(\Lm_0-\Lm_1)$. Thus,
  $\vp$ preseves
  weights of crystals. By Example \ref{Example:crystal} (ii),
  we can see that
  $\eit(t_{\lm}\ot (n))=t_{\lm}\ot\eit(n)$,
  $\fit(t_{\lm}\ot (n))=t_{\lm}\ot\fit(n)$,
  $\eit((m+n)\ot t_{-\lm})=\eit(m+n)\ot t_{-\lm}$,
  $\fit((m+n)\ot t_{-\lm})=\fit(m+n)\ot t_{-\lm}$.
  Thus we get $\vp\eit=\eit\vp$ and
  $\vp\fit=\fit\vp$. Now, it is remaining to show
  that
  $\vep_i(t_{\lm}\ot(n))=\vep_i((m+n)\ot t_{-\lm})$ and
  $\vp_i(t_{\lm}\ot(n))=\vp_i((m+n)\ot t_{-\lm})$.
  In order to show these, it is enough to notice the following
  lemma.
  \begin{lem5}
  Let $B$ be a crystal. For $b\in B$ and $\lm$, $\mu\in P$, we
  get
  \begin{eqnarray*}
  &&\vep_i(t_{\lm}\ot b\ot t_{\mu})=\vep_i(b)-\lan
  h_i,\lm\ran,\\
  &&\vp_i(t_{\lm}\ot b\ot t_{\mu})=\vp_i(b)+\lan h_i,\mu\ran.
  \end{eqnarray*}
  \end{lem5}
  This lemma is trivial from Example \ref{Example:crystal}
  (ii),
  (\ref{tensor-e}) and (\ref{tensor-f}).
  Therefore, we obtain the following formula:
  \begin{eqnarray*}
  &&\vep_1(t_{\lm}\ot(n))=m+n=\vep_1((m+n)\ot t_{-\lm}),\\
  &&\vep_0(t_{\lm}\ot(n))=-m-n=\vep_0((m+n)\ot t_{-\lm}),\\
  &&\vp_1(t_{\lm}\ot(n))=-n=\vp_1((m+n)\ot t_{-\lm}),\\
  &&\vp_0(t_{\lm}\ot(n))=n=\vp_0((m+n)\ot t_{-\lm}).
  \end{eqnarray*}
  We completed the proof. \qed

Applying
  \begin{eqnarray}
  &&T_{\lm}\ot T_{-\lm}\cong T_{-\lm}\ot T_{\lm}\cong T_0,
  \label{TT=0}\\
  &&B\ot T_0\cong T_0\ot B\cong B,
  \label{BT=TB=B}
  \end{eqnarray}
  to (\ref{TB=BT}), we get  isomorphisms:
  \begin{eqnarray}
  \vp_-:\binf & \mapright{\sim} & T_{-\lm}\ot \binf \ot T_{-\lm},
  \label{B-}\\
       (n) & \mapsto & t_{-\lm}\ot (m+n)\ot
  t_{-\lm},\nonumber\\
  \vp_+:\binf & \mapright{\sim} & T_{\lm}\ot \binf \ot T_{\lm},
  \label{B+}\\
       (n) & \mapsto & t_{\lm}\ot (n-m)\ot t_{\lm}.\nonumber
  \end{eqnarray}

  \noindent
  By applying (\ref{B-}) and (\ref{B+}) to (\ref{B=BB}), we
  get the following
  isomorphisms of crystal:
  \begin{eqnarray}
  B(-\infty) & \mapright{\sim} & T_{-\lm}\ot \binf\ot
  T_{-\lm}\ot B(-\infty),
  \label{B=TBTB-}\\
  u_{-\infty} & \mapsto & t_{-\lm}\ot(m)\ot t_{-\lm}\ot
  u_{-\infty},\nonumber\\
  B(-\infty) & \mapright{\sim} & T_{\lm}\ot \binf\ot T_{\lm}\ot
  B(-\infty),
  \label{B=TBTB+}\\
  u_{-\infty} & \mapsto & t_{\lm}\ot(-m)\ot t_{\lm}\ot
  u_{-\infty}.\nonumber
  \end{eqnarray}
  By combining (\ref{B=TBTB-}) and (\ref{B=TBTB+}), and using
  (\ref{TT=0}) and (\ref{BT=TB=B}) again, we obtain an
  isomorphism of crystal,
  \begin{eqnarray}
  T_{\lm}\ot B(-\infty) & \mapright{\sim} &  \binf\ot\binf\ot
  T_{\lm}\ot B(-\infty)
  \label{TB=BBTB}\\
  t_{\lm}\ot u_{-\infty} & \mapsto & (m)\ot (-m)\ot t_{\lm}\ot
  u_{-\infty}.
  \nonumber
  \end{eqnarray}
  Now, we set
  $$
  {\cal P}_m(-\infty):=
  \{p=(i_0,i_1,\cdot\cdot,i_k,\cdot\cdot)|i_k\in \binf{\hbox{
  and if }}|k|\gg0,
  i_{2k}=(m) {\hbox{ and }}i_{2k+1}=(-m)\}.
  $$
  By using (\ref{TB=BBTB}) repeatedly and
  arguing similarly as in \ref{5.1}, we get
  \begin{pr5}
  The following is an isomorphism of crystal;
  \begin{equation}
  T_{\lm}\ot B(-\infty)\cong {\cal P}_m(-\infty).
  \label{m-path}
  \end{equation}
  Here note that
  $t_{\lm}\ot u_{-\infty}\mapsto (m)\ot (-m)\ot (m)\ot
  (-m)\ot\cdots$.
  \end{pr5}

\vskip10pt
We set
\beq
\begin{array}{rcl}
&& {\cal P}_m :=
\left\{
p=(\cdots,i_k,i_{k+1},\cdots,i_{-1},i_0,i_1,\cdots,i_l,i_{l+1},\cdots)
|\right. i_k\in B_{\infty}  \\
&&\,\,\,\left.{\hbox{ if }}k\ll0,\,\,i_k=(0)
  {\hbox{ and if }}l\gg0,\,\, i_{2l}
=(m){\hbox{ and }}i_{2l+1}=(-m)\right\}.\nonumber
\end{array}
\label{path}
\eeq
  Now let us call an element of ${\cal P}_m$ $m$-{\it
  path} or
  simply, {\it path}.

  By applying (\ref{0-path}) and (\ref{m-path}) to Theorem
  \ref{U=BTB},
  we can easily obtain the following result:
  \newtheorem{thm5}[def5]{Theorem}
  \begin{thm5}
  \label{thm:U=path}
  For $\lm=m(\Lm_0-\Lm_1)\in P_{cl}$ $(m\in \ZZ)$,
  \begin{equation}
  \label{eqn:U=path}
  B(U'a_{\lm})\cong {\cal P}_m.
  \end{equation}
  \end{thm5}

  \vskip10pt
  \noindent
  Here note that this is an isomorphism of classical
  crystals.

\subsection{Wall and Domain}
\label{Wall and Domain}
  For this subsection, see {\it e.g.} \cite{DFJMN},
  \cite{IIJMNT}.
In the rest of this paper, we identify $B_{\infty}$ with $\ZZ$. Thus,
for $(i)\in B_{\infty}$ we denote $i$ and then for $i,j\in B_{\infty}$
we can formally consider the summation  and subtraction $i\pm j$,
and the absolute value $|i|$.

We fix an integer $m\in \ZZ$ and
  let $p\in {\cal P}_m$ be a $m$-path.
  \begin{def5}
  \label{wall}
\begin{enumerate}
  \item
  A path $p=(\cdots, i_{k-1},i_k,\cdots)$ has $l$ walls in the
  position $k$ $($$l\in \ZZ_{> 0}$, $k\in \ZZ$$)$, if
  $|i_{k-1}+i_k|=l$.
  \item
  Suppose that there are walls in position $k$. The {\it type
  of walls} in position $k$ is $+$ $($resp. $-$$)$ if
  $i_{k-1}+i_k>0$
  $($resp. $i_{k-1}+i_k<0$$)$.
  \end{enumerate}
  \end{def5}

  \noindent
  We also define a function $n:{\cal P}_m\longrightarrow \ZZ_{\geq0}$  by
  $$
  n(p)=\sum_{k\in\ZZ}|i_{k-1}+i_k|
  $$
  and we  call this the total number of walls in $p$.

  Here note that for any $p\in {\cal P}_m$, $n(p)$ is
  finite
  by the definition of ${\cal P}_m$.
  \begin{def5}
  \label{def-domain}
  A segment $d= i_j,i_{j+1},\cdots,i_l$ in $p\in {\cal P}_m$ is a
  {\it finite domain} with length $l-j+1$ in $p$ if there are
  walls in the position $j$ and $l+1$ and there is no wall in
  positions $j+1,j+2,\cdots,l$. We denote $l(d):=l-j+1$ for
the length of domain $d$.
  \end{def5}

  \noindent
  {\sl Remark.}
  \begin{enumerate}
  \item
  In this definition, we can consider a domain with length 0.
  This occurs in the following case. If there are more than
  two walls in the same position,  there is a domain with
  length 0 between a pair of
  neighboring two walls.
  \item
  By the definition of $\mpath$, we know that any path has two
  infinite sequences in
  the forms $\cdots0,0,0,0$ and $\pm m, \mp m, \pm m,\cdots$.
  We call these {\it infinite domains}.
  \item
  By the definition of finite domain,
any finite domain with positive length is in
  the following form;
  \beq
  k,-k,k,-k,\cdots, \pm k,\mp k.
  \eeq
  \end{enumerate}

  \newtheorem{ex5}[def5]{Example}
  \begin{ex5}
  For $p=(\cdots,0,0,1,-1,3,-3,3\cdots)$,
we visualize walls and domains:
\beq
\cdots 00|1-1||3-33\cdots.
\label{vis}
\eeq
In (\ref{vis}),
we know that there are three walls, two finite domains $:$
$1-1$ and a zero-length domain and two infinite domains $:$
$\cdots 00$ and $3-33\cdots$.
  \end{ex5}

  \vskip5pt
  \noindent
  Now, for $n\in \ZZ_{\geq 0}$ we set
  $$
  \mpath(n):=\{p\in \mpath|n(p)=n\}.
  $$
  It is trivial that $\mpath=\oplus_{n\geq0}\mpath(n)$.
  By simple calculations, we get
 \begin{pr5}
\label{empty-condition}
\begin{enumerate}
  \item
  If $m$ is odd $($resp. even$)$, then $\mpath(2n)=\emptyset $
  $($resp. $\mpath(2n-1)=\emptyset$$)$.
  \item
  If $n<|m|$, then $\mpath(n)=\emptyset$.
  \end{enumerate}
\end{pr5}

  We shall see the stability of $\mpath(n)$ by the actions of
  $\eit$ and $\fit$.
  \begin{pr5}
  \label{pmn}
  For a path $p\in \mpath(n)$, suppose that $\fit p\ne 0$
  $($resp. $\eit p\ne 0$$)$, then we have $n(\fit p)=n(p)$ $($resp.
  $n(\eit p)=n(p)$$)$.
  \end{pr5}

  \noindent
  {\sl Proof. }\,\,
  For a path $p=(\cdots, i_k,i_{k+1}\cdots)$ and $i=0,1$, we
  set
  \begin{equation}
  a^{(i)}_k=\sum_{j<k}\vp_i(i_j)-\vep_i(i_{j+1}).
\label{a_k}
  \end{equation}

\nd
{\sl Remark.}\,
If $k\ll0$ then $i_k=0$,  thus we have
  $a^{(i)}_k=0$ for $k\ll0$ and by the fact that $\vp_i(\pm
  m)=\vep_i(\mp m)$ we have $a^{(i)}_k=a^{(i)}_{k+1}$ for
  $k\gg0$.

In order to prove the proposition, it is necessary to see
  the following lemma.
\begin{lem5}
\label{action}
  \begin{enumerate}
  \item
  For a path $p=(\cdots, i_k,i_{k+1}\cdots)$, if
  there exists
  $k\in \ZZ$ such that
  \begin{equation}
  a^{(i)}_{\nu}\geq a^{(i)}_k \,\,(\nu<k){\hbox{ and }}
  a^{(i)}_{\nu} > a^{(i)}_k \,\,(\nu>k),
\label{ai:fi}
\end{equation}
  then
\begin{equation}
\fit p=(\cdots,i_{k-1}, \fit(i_k),i_{k+1}\cdots),
\label{fip}
\end{equation}
  otherwise $\fit p=0$.
\item
  For a path $p=(\cdots, i_k,i_{k+1}\cdots)\in \mpath$, if
  there exists
  $k\in \ZZ$ such that
\begin{equation}
  a^{(i)}_{\nu}>a^{(i)}_k \,\,(\nu<k){\hbox{ and }}
  a^{(i)}_{\nu} \geq a^{(i)}_k \,\,(\nu>k),
\end{equation}
  then
\begin{equation}
\eit p=(\cdots,i_{k-1}, \eit(i_k),i_{k+1}\cdots),
\end{equation}
  otherwise $\eit p=0$.
 \end{enumerate}
  \end{lem5}

  \noindent
  {\sl Proof of Lemma \ref{action}}\,\,
Since the proof of (ii) is similar to the one of (i),
we shall show only (i).
 For any $p=(\cdots,i_{k},i_{k+1},\cdots)\in \mpath(n)$
 there exist $j,l\in\ZZ_{>0}$ such that $i_k=0$
if $k\leq -j$ and $i_{2k}=m$ and
$i_{2k+1}=-m$ if $k\geq l$. Then $p$ is identified with
\beq
u_{\infty}\ot i_{-j}\ot i_{-j+1}\ot \cdots
\ot i_{2l}\ot i_{2l+1}\ot t_{\lm}\ot u_{-\infty}.
\label{infinf}
\eeq
For the vector (\ref{infinf}),
we set $A^{(i)}_{-j}:=0$ and
\beq
A^{(i)}_k=\sum_{-j\leq p<k}\vp_i(i_p)-\vep_i(i_{p+1}),
\,\,\,\,\,k=-j+1,\cdots,2l.
\label{AK}
\eeq
Since $\vp_i(u_{\infty})
=\vep_i(u_{\infty})=\vp_i(i_{-j})=\vep_i(i_{-j})=0$,
we get
\beq
A^{(i)}_k=a^{(i)}_k \,\,{\hbox{ for $k=-j,\cdots,2l$}}.
\label{A=a}
\eeq
By the previous remark,
if there exists $k$ which satisfies (\ref{ai:fi}),
$-j\leq k\leq 2l$.
Therefore, by Proposition 2.1.1 (i) in \cite{KN}
we obtain the desired result.\qed

Now, let us show Proposition \ref{pmn} (i).
We shall consider $i=1$ case.
Suppose that for
$p=(\cdots,i_{k-1},i_k,i_{k+1}\cdots)$ we have
$\til f_1\,p=(\cdots,i_{k-1},\til f_1(i_k),i_{k+1}\cdots)$.
We know that $\til f_1(i_k)=i_k+1$. Thus, we get
$$
\til f_1\,p=(\cdots,i_{k-1},i_k+1,i_{k+1}\cdots).
$$
By Lemma \ref{action}, we have
$$
a^{(1)}_{k-1}\geq a^{(1)}_{k}{\hbox{ and }}a^{(1)}_{k+1}>a^{(1)}_k.
$$
By using this, we obtain,
\begin{eqnarray}
0\leq a^{(1)}_{k-1}-a^{(1)}_k
& = & -(\vp_1(i_{k-1})-\vep_1(i_k)) =i_{k-1}+i_k,
\label{pmn-a}\\
0< a^{(1)}_{k+1}-a^{(1)}_k
& = & \vp_1(i_k)-\vep_1(i_{k+1}) =-i_{k}-i_{k+1}.
 \label{pmn-b}
\end{eqnarray}
By (\ref{pmn-a}) and (\ref{pmn-b}), we get
$$
i_{k-1}+i_k\geq0\,\,{\hbox{ and }}\,\,i_{k}+i_{k+1}< 0.
$$
Therefore,
\begin{eqnarray*}
|i_{k-1}+\til f_1(i_k)| & = & |i_{k-1}+i_k+1|=|i_{k-1}+i_k|+1, \\
|\til f_1(i_k)+i_{k+1}| & = & |i_k+i_{k+1}+1|=|i_k+i_{k+1}|-1.
\end{eqnarray*}
Then $n(\til f_1\, p)=n(p)$.
By arguing similarly we can prove for $i=0$.
Now, we have completed the proof of Proposition \ref{pmn}.
\hfill\framebox[2mm]{}

%%%%%%%% section 6 %%%%%%%%%
\section{Path-Spin Correspondence}
\setcounter{equation}{0}
\renewcommand{\theequation}{\thesection.\arabic{equation}}

The purpose of this section is to give a
strict morphism of $P_{cl}$-weighted crystals
$\mpath(n)\longrightarrow B^{\ot n}$. ({\it cf.}\,  2.2)

Now,  we shall define a map from $\mpath(n)$ to $B^{\ot n}$ as follows:
For $p\in\mpath(n)$, let $(\io_1,\io_2,\cdots,\io_n)$ be
the sequence of wall types (ordered from the left to the right).
The map $\psi:\mpath(n)\longrightarrow B^{\ot n}$ is given by
\beq
\psi(p)=(-\io_1)\ot (-\io_2)\ot\cdots\ot(-\io_n),
\label{V}
\eeq
for any $p\in \mpath(n)$.
\newtheorem{thm6}{Theorem}[section]
\begin{thm6}
\label{path-spin}
The map $\psi$ is
a strict morphism of $P_{cl}$-weighted crystals from
$\mpath(n)$ to $B^{\ot n}$.

\end{thm6}

\noindent
{\sl Proof.}\,\,
In order to prove the theorem, we shall see that
$\psi$ satisfies
\beqn
 wt(p)&=&wt(\psi(p)),
\label{wt-psi} \\
\vep_i(p)&=&\vep_i(\psi(p)),\,\,\vp_i(p)=\vp_i(\psi(p)),
\label{vep-vp-psi}\\
\fit\psi(p) & = & \psi(\fit \, p),
\label{fipsi=psifi} \\
\eit\psi(p) & = & \psi(\eit \, p),
\label{eipsi=psiei}
\eeqn
for any $p\in \mpath(n)$ and $i=0,\,1$.

\vskip7pt
\noindent

A $m$-path $g=(g_k)_{k\in \ZZ}$
satisfying $g_k=0$ for $k<0$, $g_{2k}=m$
and $g_{2k+1}=-m$
for $k\geq0$ is called $m$-{\it ground-state path}.
$g=(g_k)_{k\in \ZZ}$ just corresponds to
 $u_{\infty}\ot t_{\lm}\ot u_{-\infty}$
in $B(\infty)\ot T_{\lm}\ot B(-\infty)$.
Then $wt(g)=m(\Lm_0-\Lm_1)$.
Therefore, for $p=(i_k)_{k\in\ZZ}$
the following formula is obtained easily;
\beq
\begin{array}{rcl}
wt(p) & = & m(\Lm_0-\Lm_1)+\sum_{k\in\ZZ}(wt(i_k)-wt(g_k))
\label{path-wt-cl} \\
   & = & (m+2 \sum_{k\in \ZZ}(i_k-g_k))(\Lm_0-\Lm_1).\nonumber
\end{array}
\eeq
By the definition of path,
we know that the summation in (\ref{path-wt-cl})
is finite. Therefore, by the fact $g_{k-1}+g_k=0$ $(k\ne 0)$ and
$g_{-1}+g_0=m$,
\beqn
wt(p)
& = & (m+\sum_{k\in \ZZ}(i_{k-1}+i_{k}-g_{k-1}-g_{k})(\Lm_0-\Lm_1)
\nonumber \\
& = & (\sharp\{(+){\rm walls}\}-\sharp\{(-){\rm walls}\})(\Lm_0-\Lm_1)
\nonumber \\
      & = & wt(\psi(p)).
\nonumber
\eeqn
Here note that $wt(\pm)=\pm(\Lm_1-\Lm_0)$.  Now we get (\ref{wt-psi}).

\vskip6pt
\noindent
Let us show (\ref{vep-vp-psi}).
For $p=(\cdots,i_k,i_{k+1},\cdots)\in \mpath(n)$,
let $a$ and $b$ be sufficiently large integers such that
$i_{-j}=0$, $i_{2k}=m$ and $i_{2k+1}=-m$
for any $j>a$ and $k>b$. Therefore, since $p$ is identified with
$u_{\infty}\ot i_{-j}\ot
\cdots\ot i_{2k}\ot i_{2k+1}\ot t_{\lm}\ot u_{-\infty}$ and
$\vp_i(u_{-\infty})=\vep_i(u_{-\infty})=0$,
by (\ref{tensor-vp}) we have
\beq
\vp_i(p)=
\vp_i(u_{\infty}\ot i_{-j}\ot
\cdots\ot i_{2k}\ot i_{2k+1}\ot t_{\lm}),\nonumber
\eeq
for $j>a$ and $k>b$.
By the formula $\vp_i(t_{\lm})=-\infty$,
Proposition 2.1.1 (0) in \cite{KN}
and (\ref{tensor-vp}), we get
\beq
\vp_i(p)
=\lan h_i,\lm\ran+\vp_i(i_{2k+1})
 +\mathop{{\rm max}}_{-j\leq p\leq 2k+1}
(a^{(i)}_{2k+1}-a^{(i)}_p ).
\label{vpip}
\eeq
We shall consider $i=1$ case.
Then (\ref{vpip}) can be written explicitly as follows:
\beq
\vp_1(p) =
\mathop{{\rm max}}_{-j\leq p\leq 2k+1}
 (-\sum_{p< s\leq 2k+1}i_{s-1}+i_{s}),
\label{phi1max}
\eeq
by using $\vp_1(i_{2k+1})=-i_{2k+1}= m = - \lan h_1,\lm \ran $.

Let $k_1,k_2,\cdots,k_s$ ($s\leq n$) be the sequence of
positions of walls in $p$ such that $k_j<k_{j+1}$
and there is no wall in
$j\ne k_1,\cdots,k_s$.
Here note that since more than
one walls can occupy the same position, $s\leq n$.
Let $c_i$ be the position of $i$-th wall.
For $j=1,2,\cdots\ s$, we set
$$
N^{\pm}_j  :=
\sharp\left\{\io_r=\pm\,|\,c_r \in \{k_j,\cdots,k_s\}\right\}.
$$
Since $i_{c-1}+i_c=0$ if $c\not\subset\{k_1,\cdots,k_r\}$,
The formula (\ref{phi1max}) can be written as follows;
\beq
\begin{array}{rcl}
\vp_1(p) & = &
\mathop{{\rm max}^*}_{1\leq j\leq s}
\left\{-\sum_{l=j}^si_{k_l-1}+i_{k_l}\right\},
 \\
& = &\mathop{{\rm max }^*}_{1\leq j\leq s}\{N^-_j-N^+_j\},
\nonumber
\end{array}
\label{walls}
\eeq
where ${\rm max}^*\{z_1,\cdots,z_n\}
:={\rm max}\{z_1,\cdots,z_n,0\}\geq 0$.
Note that if there is no $(-)$ wall in $p$, $\vp_1(p)=0$ and
$\vp_1(\psi(p))=\vp_1((-)^{\ot n})=0$.
Then we may assume that there exists $(-)$ wall in $p$.

\noindent
We shall investigate $\vp_1(\psi(p))$.
By Proposition 2.1.1 (0) in \cite{KN}, we can get the following,
\beqn
\vp_1(\psi(p)) & = & \vp_1((-\io_1)\ot \cdots\ot(-\io_n))
\nonumber \\
& = &
\mathop{{\rm max}}_{1\leq j\leq n}
\left\{\sum_{j\leq k\leq n}\vp_1(-\io_k)
-\sum_{j<k\leq n}\vep_1(-\io_k)\right\}.
\label{epep}
\eeqn
Here note that $\vp_1(+)=1=\vep_1(-)$ and
$\vp_1(-)=0=\vep_1(+)$ by (\ref{act}),
$\sum_{j\leq k\leq n}\vp_1(-\io_k)
=\sharp\{\io_k=-\,;\, j\leq k\leq n\}$ and
$\sum_{j< k\leq n}\vep_1(-\io_k)
=\sharp\{\io_k=+\,;\, j< k\leq n\}$.
Then we can rewrite (\ref{epep}) as follows;
\beq
\vp_i(\psi(p))=\mathop{{\rm max}}_{1\leq j\leq n}
\{\sharp\{\io_k=-\,;\, j\leq k\leq n\}-
\sharp\{\io_k=+\,;\, j< k\leq n\}\}.
\label{+--}
\eeq
Therefore, if $t$ ($1\leq t \leq n$)
gives the maximum in (\ref{+--}), there are two cases
\begin{enumerate}
\item
$\io_t=-$ and $\io_{t-1}=+$ $(t>1)$.
\item
$t=1$ and $\io_1=-$.
\end{enumerate}
Since in both cases $\vep_1(-\io_t)=\vep_1(+)=0$ and
$\sum_{j\leq k\leq n}\vep_1(-\io_k)\geq
 \sum_{j< k\leq n}\vep_1(-\io_k)$,
we can rewrite
(\ref{+--}) to
\beq
\vp_i(\psi(p))=\mathop{{\rm max}}_{1\leq j\leq n}
\{\sharp\{\io_k=-\,;\, j\leq k\leq n\}-
\sharp\{\io_k=+\,;\, j\leq k\leq n\}\}.
\label{+---}
\eeq
Since we have the following by the definition of $N^{\pm}_j$,
$$
\left\{N^-_j-N^+_j\right\}_{1\leq j\leq s}\subset
\{\sharp\{\io_k=-\,;\, j\leq k\leq n\}-
\sharp\{\io_k=+\,;\, j\leq k\leq n\}\}_{1\leq j\leq n},
$$
by (\ref{walls}) and (\ref{+---})  we get $\vp_1(p)\leq\vp_i(\psi(p))$.
We set
\beq
S:=\{s \,|\,
\hspace{-8pt}
\begin{array}{ll}
&1\leq s\leq  n,{\hbox{ $s$-th wall in $p$ is a $(-)$ wall and }}\\
&{\hbox{the left-most wall among walls in the same position}}
\end{array}\}.
\label{SSS}
\eeq
The cases (i) and (ii) as above
mean that if $t$ gives the maximum of (\ref{+---}), $t\in S$.
Here note that if $s\in S$,
$$
N^{\pm}_s=\sharp\{\io_k=\pm;s\leq k\leq n\}.
$$
Therefore, we get $\vp_1(p)\geq\vp_1(\psi(p))$.
Now, we have $\vp_1(p)=\vp_1(\psi(p))$.
As for $\vp_0$-case and $\vep_i$-case arguing similarly,
we obtain (\ref{vep-vp-psi}).

Let us show (\ref{fipsi=psifi}) for $i=1$.
For $p=(\cdots,i_{j-1},i_j,i_{j+1}\cdots)$ we assume that
there exists $k$ satisfying (\ref{ai:fi}) for $i=1$, {\it i.e.}
$\til f_1\,p=(\cdots,i_{k-1},\til f_1(i_{k}),i_{k+1},\cdots)$.

\noindent
We know that $a^{(1)}_k$ is given by
$a^{(1)}_k=-\sum_{j<k}i_j+i_{j+1}$.
Since $k$ satisfies (\ref{ai:fi}) for $i=1$, we have
$ a^{(1)}_k< a^{(1)}_{k+1}$ and  $a^{(1)}_{k-1}\geq a^{(1)}_k$.
Then we get
\beqn
&& i_k+i_{k+1}=a^{(1)}_k-a^{(1)}_{k+1}<0,
\label{ik+ik+1<0} \\
&& i_{k-1}+i_k=a^{(1)}_{k-1}-a^{(1)}_k\geq 0.
\label{ik-1+ik>0}
\eeqn
By (\ref{ik+ik+1<0}), we know that
there exist $(-)$ walls in position $k+1$.
By this and the fact
\beq
\til f_1(i_k)=i_k+1,
\label{f1ik=ik+1}
\eeq
we get
\beq
|\til f_1(i_k)+i_{k+1}|=|i_k+i_{k+1}|-1.
\label{-wallheru}
\eeq
By (\ref{ik-1+ik>0})  and (\ref{f1ik=ik+1}),
\beq
|i_{k-1}+\til f_1(i_k)|=|i_{k-1}+i_k|+1,
\label{+wallfueru}
\eeq

\nd
Let $j$-th wall in $p$ be the left-most wall among walls
in position $k+1$ (the existence is
guaranteed by (\ref{ik+ik+1<0})).
(\ref{ik+ik+1<0}), (\ref{-wallheru}) and (\ref{+wallfueru})
imply that the $j$-th wall and
other walls in position $k+1$
are  $(-)$ walls and
the $j$-th wall is changed by the action of
$\til f_1$ to $(+)$ wall in position $k$.
Then $j$ belongs to $S$.

That is, let $(\io_1,\io_2,\cdots,\io_n)$  and
$(\io'_1,\io'_2,\cdots,\io'_n)$ be
the sequences of wall types of $p$ and $\til f_1\, p$
respectively, we have
\beq
(\io_1,\cdots,\stackrel{j}{-},\cdots,\io_n)\mapright{\til f_1}
(\io_1,\cdots,\stackrel{j}{+},\cdots,\io_n)
=(\io'_1,\cdots,\io'_n).
\label{-f1+}
\eeq
By (\ref{-f1+}), we know that
\beq
\psi(p)= (-\io_1)\ot
\cdots \ot \stackrel{j}{+}\ot\cdots\ot(-\io_n),\,\,\,
\psi(\til f_1p)=(-\io_1)\ot
\cdots \ot \stackrel{j}{-}\ot\cdots\ot(-\io_n).
\label{psippsif1p}
\eeq
By (\ref{psippsif1p}), it is sufficient to show the following,
\beq
\til f_1\left((-\io_1)\ot
\cdots\ot \stackrel{j}{+}\ot\cdots\ot(-\io_n)\right)=
(-\io_1)\ot\cdots\ot \stackrel{j}{-}\ot\cdots\ot(-\io_n).
\label{BBB}
\eeq
For $p$ with $\psi(p)=(-\io_1)\ot\cdots\ot(-\io_n)$
we shall define the function $\overline{a}_k$ as follows;
(this just coincides with $a_k$
in Proposition 2.1.1 (0) in \cite{KN} up to the first term.).
\beq
\overline{a}_k:=-\vep_1(-\io_1)+
\sum_{1\leq l< k}\vp_1(-\io_l)-\vep_1(-\io_{l+1}).
\label{abar}
\eeq
It is easy to translate (\ref{abar}) to the following form
by (\ref{act});
\beq
\overline{a}_k=\sharp\{\io_l=-\,|\, 1\leq l<k\}
-\sharp\{\io_l=+\,|\, 1\leq l\leq k\}.
\label{no-+}
\eeq
By Proposition 2.1.1 (i) in \cite{KN}, we know that
if there exists $j$ satisfying
\beq
\overline{a}_{\nu}\geq \overline{a}_j
 \,\,{\hbox{ for }}\,\,\nu<j\,\,{\hbox{ and }}
\overline{a}_{\nu}   > \overline{a}_j
\,\,{\hbox{ for }}\,\, j<\nu,
\label{anuaj}
\eeq
we have
\beq
\til f_1\left((-\io_1)\ot\cdots\ot(-\io_j)\ot\cdots\ot(-\io_n)\right)
=(-\io_1)\ot\cdots\ot\,\til f_1(-\io_j)\ot\cdots\ot(-\io_n).
\label{fi-ioj}
\eeq
Then we shall show that $j$ as in (\ref{-f1+})
and (\ref{psippsif1p}) satisfies (\ref{anuaj}).
Since the position of the $j$-th wall is $k+1$ and
there is no $(+)$ wall in position $k+1$ by the argument above,
we get
\beqn
\sharp\{\io_l=-\,|\, 1\leq l<j\}
& = & \hspace{-8pt}\sum_{\stackrel{\scriptstyle r\leq k}{i_{r-1}+i_r<0}}
 |i_{r-1}+i_r|
  =  -\sum_{\stackrel{\scriptstyle r\leq k}{i_{r-1}+i_r<0}}
 i_{r-1}+i_r,
\label{io-} \\
\sharp\{\io_l=+\,|\, 1\leq l\leq j\}
& = & \hspace{-8pt}\sum_{\stackrel{\scriptstyle r\leq k+1}{i_{r-1}+i_r>0}}
 |i_{r-1}+i_r|
  =   \sum_{\stackrel{\scriptstyle r\leq k}{i_{r-1}+i_r>0}}
 i_{r-1}+i_r.
\label{io+}
\eeqn
The following is obtained by (\ref{no-+}), (\ref{io-}) and (\ref{io+}),
\beq
\overline a_j =  -\sum_{r\leq k}i_{r-1}+i_r
            = \sum_{r< k}\vp_1(i_r)-\vep_1(i_{r+1})  = a^{(1)}_k,
\label{aj=ak}
\eeq
By the form of (\ref{no-+}),
we know that $j$ satisfying (\ref{anuaj}) belongs to $S$.
Therefore, in order to show that $j$ satisfies (\ref{anuaj})
it is enough to show that $j$ satisfies
\beq
\overline{a}_{\nu}\geq \overline{a}_j
\,\,{\hbox{ for }}\,\,\nu<j\,\,(j,\nu\in S)\,{\hbox{ and }}
\overline{a}_{\nu}   > \overline{a}_j
\,\,{\hbox{ for }}\,\, j<\nu,\,\,(j,\nu\in S).
\label{anuajs}
\eeq
By the same argument as for obtaining (\ref{aj=ak}),
we can see that
for any $\nu<j$ (resp. $\nu>j$) ($\nu,j\in S$)
there exists $t$ such that
\beq
t<k\,\,({\rm resp. }\,\,t>k) \,\,{\hbox{ and }}\,\,
\overline a_{\nu}=a^{(1)}_t.
\label{as=at}
\eeq
By (\ref{ai:fi}) for $i=1$,
(\ref{aj=ak}) and (\ref{as=at}), we get that
$j$ satisfies (\ref{anuajs}) and then (\ref{anuaj}).
Now, we get (\ref{fi-ioj}).

Next, we shall show that if $\til f_1p=0$, $\til f_1\psi(p)=0$.
We assume that $\til f_1p=0$ and set
$\xi={\rm max}\{\nu\,|\,i_{\nu-1}+i_{\nu}\ne 0\}$.
By Lemma \ref{action} we know that
$\xi$ satisfies
\beq
a^{(1)}_{\nu}\geq a^{(1)}_{\xi}\,\,{\hbox{ for $\nu>\xi$
and $a^{(1)}_{\xi}= a^{(1)}_{\nu}$ for $\xi>\nu$.}}
\label{anuuu}
\eeq
Now, we set $F:=a^{(1)}_{\xi}$.
Let us assume that $a^{(1)}_{\xi-1}=a^{(1)}_{\xi}$.
Then we have
$$
0=a^{(1)}_{\xi-1}-a^{(1)}_{\xi}
=-\vp_1(i_{\xi-1})+\vep_1(i_{\xi})=i_{\xi-1}+i_{\xi}.
$$
This contradicts  the definition of $\xi$. Thus, we get
$a^{(1)}_{\xi-1}>a^{(1)}_{\xi}$ and then $i_{\xi-1}+i_{\xi}>0$.
Furthermore, by the fact that $i_{\xi-1}+i_{\xi}>0$ we have
\beq
\io_n=+.
\label{ion=+}
\eeq
Here note that
\beq
a^{(1)}_{\xi}=F=\overline a_n.
\label{amu=an}
\eeq
It is sufficient to show that
\beq
\til f_1\psi(p)=\til f_1\left((-\io_1)\ot\cdots\ot(-\io_n)\right)
=(-\io_n)\ot\cdots\ot\til f_1(-\io_n),
\label{f1psip}
\eeq
since we know that (\ref{ion=+}) and $\til f_1(-)=0$.
Since we know that if there is $j$ with (\ref{anuaj}),
(\ref{fi-ioj}) holds.
in order to show (\ref{f1psip}),
we shall prove
\beq
\overline a_{\nu}\geq \overline a_n,\,\,{\hbox{ for $\nu< n$.}}
\label{anuan}
\eeq

\noindent
We assume that there exists $j$ such that $j\ne n$
and satisfies (\ref{anuaj}).
Let $t$ be the position of the $j$-th wall.
It is easy to see that  $\io_{j}=-$ by (\ref{no-+}).
Thus, by (\ref{ion=+}), we have
\beq
t<\xi.
\label{t<mu}
\eeq

By simliar argument to
the one for obtaining (\ref{as=at}),
we get $\overline a_j=a^{(1)}_{t}$.
Therefore, by (\ref{anuuu}) and (\ref{amu=an}) we have
\beq
\overline a_j\geq F=\overline a_n.
\label{aj>F}
\eeq
This contradicts  the definition of $j$ satisfying (\ref{anuaj}).
Now we get (\ref{f1psip}) and
then $\til f_1\psi(p)=0$ if $\til f_1p=0$.

By arguing similarly, we obtain
$\til f_0\psi(p)=\psi(\til f_0p)$ and
$\til e_i\psi(p)=\psi(\til e_ip)$.
Then, we have completed the proof of Theorem \ref{path-spin}. \qed

%%%%%%%%%%% section 7 %%%%%%%%%%%
\section{Classification of Path}
\setcounter{equation}{0}
\renewcommand{\theequation}{\thesection.\arabic{equation}}

In this section, we shall describe every connected component in $B(\util)$.
The notations in this section follow the previous sections.

\subsection{Domain type and Domain parameter}
For a path $p\in \mpath(n)$ $(n\geq 0,\,\,m\in \ZZ)$,
let $d_0,d_1,\cdots,d_{n-1},d_n$ be the sequence of domains in $p$.
The domains $d_0$ and $d_n$ are infinite domains.
Here note that as stated in Remark (iii) in \ref{Wall and Domain},
\beq
{\hbox{every domain with non-zero length is in the form:\,\, }}
\cdots,k,-k,k,-k,\cdots,
\label{l-l}
\eeq
where $k$ is an integer.

\newtheorem{def7}{Definition}[section]
\begin{def7}
For a domain $d_j$ with non-zero length,
fixing some entry $i_{\nu}$ in $d_j$ and its position $\nu$,
the {\it domain type} $t(d_j)$ of $d_j$ is given by
\beq
t(d_j):=(-1)^{\nu}i_{\nu}.\,\,
\eeq
\end{def7}
{\sl Remark.\,}
\begin{enumerate}
\item
By (\ref{l-l}),
this definition is well-defined, {\it i.e.},
a domain type is uniquely determined.
\item
Domain type of domain $d_0$ is always 0 and
one of domain $d_n$ is always $m$ by the definition of $\mpath(n)$.
\end{enumerate}

The following lemma is trivial.
\newtheorem{lem7}[def7]{Lemma}
\begin{lem7}
\label{ab}
For a path $p$ let $i_{k-1}$ and $i_k$ be entries
in $p$ with $|i_{k-1}+i_k|\ne 0$ and
let $d_j$ and $d_l$ ($j<l$) be domains
including $i_{k-1}$ and $i_k$ respectively.
Then we have
$$
|t(d_l)-t(d_j)|-1={\hbox{the number of domains
with zero-length between $d_j$ and $d_l$.}}
$$
\end{lem7}
By this lemma, the following definition is well-defined.
\begin{def7}
\label{zero-dom}
Let $d_r$ be the $i$-th zero-length domain
between $d_j$ and $d_l$
as in Lemma \ref{ab}. Domain type $t(d_r)$ is given by:
$t(d_r)=t(d_j)+i$ if $t(d_j)<t(d_l)$ and
$t(d_r)=t(d_j)-i$ if $t(d_j)>t(d_l)$.
\end{def7}
\newtheorem{ex7}[def7]{Example}
\begin{ex7}
\label{1212}
For $p=(\stackrel{{\rm position}}{\cdots\cdot\cdot},
\stackrel{-2}{0},\stackrel{-1}{0},\stackrel{0}{2},
\stackrel{1}{-1},\stackrel{2}{3},\stackrel{3}{-3},
\stackrel{4}{3},\cdots)$,
we shall visualize walls and domains as follows:
$$
\cdots \stackrel{d_0}{00}\stackrel{d_1}{|\,\,|}\stackrel{d_2\,\,\,d_3}{2|-1}
\stackrel{d_4}{|\,\,|}\stackrel{d_5}{3  -33}\cdots.
$$
We know that there are five walls and four finite domains in $p$.
Let $d_1$, $d_2$, $d_3$ and $d_4$ be the four finite domains.
The domains $d_1$ and $d_4$ are zero-length domains. The domain type of these
four domains are $1,2,1,2$ respectively.
Of course, the domain type of the leftmost infinite domain is 0 and the one
of the rightmost infinite domain is 3.
\end{ex7}

\noindent
{\it Remark.\,}
Note that for any path $p\in \mpath(n)$ and $j=0,1,\cdots, n-1$
\beq
|t(d_{j+1})-t(d_j)|=1.
\label{tj+1-tj=1}
\eeq
\begin{def7}
A sequence of integers $t_1,t_2,\cdots,t_{n-1}$ is
in $m$-domain configuration
if $|t_{j}-t_{j-1}|=1$ for $j=1,\cdots,n$,
where $m\in\ZZ$, $t_0=0$ and $t_n=m$.
\end{def7}
The following lemma is trivial.
\begin{lem7}
\label{n-m}
There exists a sequence $t_1,\cdots, t_{n-1}$
in $m$-domain configuration if and only if
$n-|m|\in 2\ZZ_{\geq 0}$.
\end{lem7}
By the above remark, we get
\begin{lem7}
A sequence of domain types for any
path in $\mpath$ is in $m$-domain configuration.
\end{lem7}

\begin{def7}
\label{def:domaintype}
\begin{enumerate}
\item
Let $\vec t=(t_1,t_2,\cdots,t_{n-1})$ be
in a $m$-domain configuration,
\begin{enumerate}
\item
$t_j$ is regular in $\vec t$
if $t_{j-1}-t_j=t_j-t_{j+1}$.
\item
$t_j$ is up $($resp. down$)$-regular in $\vec t$
if $t_j$ is regular in $\vec t$ and $t_{j-1}<t_j<t_{j+1}$
 $($resp. $t_{j-1}>t_j>t_{j+1}$$)$.
\item
$t_j$ is critical in $\vec t$
if $t_{j-1}-t_j=-t_j+t_{j+1}$.
\item
$t_j$ is maximal $($resp. minimal$)$ in $\vec t$
if $t_j$ is critical in $\vec t$ and $t_{j-1}+1=t_j=t_{j+1}+1$
$($resp. $t_{j-1}-1=t_j=t_{j+1}-1$$)$.
\end{enumerate}
Here $t_0=0$ and $t_n=m$.
\item
For a path $p\in \mpath(n)$, let
$d_1,\cdots,d_{n-1}$ be its finite domains
and $\vec t(\vec d)=(t(d_1),\cdots,t(d_{n-1}))$
be the sequence of their domain types.
\begin{enumerate}
\item
$d_j$ is a  regular domain
if $t(d_j)$ is regular in $\vec t(\vec d)$.
\item
$d_j$ is  {\sl up-regular} $($resp.
{\sl down-regular}$)$ if $t(d_j)$ is up-regular in
$\vec t(\vec d)$,
in particular,
$d_0$ is  {\sl up} $($resp. {\sl down}$)$
if $t(d_0)<t(d_1)$ $($resp. $t(d_0)>t(d_1)$ and
$d_n$ is  {\sl up} $($resp. {\sl down}$)$
if $t(d_{n-1})<t(d_n)$ $($resp. $t(d_{n-1})>t(d_n)$$)$.
\item
$d_j$ is a  critical domain
if $t(d_j)$ is critical in $\vec t(\vec d)$.
\item
$d_j$ is {\sl maximal} (resp. {\sl minimal})
if  $t(d_j)$ is maximal  {$($resp. minimal$)$ in $\vec t(\vec d)$.}
\end{enumerate}
\end{enumerate}
\end{def7}
{\sl Remark.}\,
\begin{enumerate}
\item
By Definition \ref{zero-dom}, any zero-length domain is a regular domain.
\item
If $t_1\cdots, t_{n-1}$ is in $m$-domain configuration,
any $t_j$ ($j=1,\cdots,n-1$) is classified
in Definition \ref{def:domaintype} (i) (b)(d) and
then any domain is classified in
Definition \ref{def:domaintype} (ii) (b)(d).
\end{enumerate}
\begin{ex7}

In Example \ref{1212}, the infinite domains $d_0$ and $d_5$ are up.
$d_1$ and $d_4$ are up-regular,
$d_2$ is maximal and $d_3$ is minimal.
\end{ex7}

\begin{def7}
For $p\in \mpath(n)$,
let $d_1,d_2,\cdots,d_{n-1}$ be its finite domains and
$l(d_1),l(d_2),\cdots,l(d_{n-1})$ be their lengths.
Domain parameter $c(d_j)$ is given by
\beqnn
&& {\hbox{if $d_j$ is a regular domain, }}\,c(d_j)
:=[[\frac{l(d_j)}{2}]],\\
&& {\hbox{if $d_j$ is a critical domain, }}\,c(d_j)
:=[[\frac{l(d_j)+1}{2}]],
\eeqnn
where $[[n]]=$the maximum integer
which is less than or equal to $n$.
\end{def7}

Let $\vec t=(t_1,t_2,\cdots,t_{n-1})$ be
in a $m$-domain configuration
and $\vec c=(c_1,c_2,\cdots,c_{n-1})$ be
a sequence of non-negative integers.
For $\vec t$ and $\vec c$, we set
\beqn
\lefteqn{\mpath(n;\vec t;\vec c)}\nonumber\\
&:=&\left\{p\in \mpath(n)\,\left|\right.
\hspace{-26pt}
{\begin{array}{rcl}
&&t(d_j)=t_j{\hbox{ and }}c(d_j)=c_j
\,{\hbox{ for any }}j=1,2,\cdots,n-1,
\nonumber\\
&&{\hbox{ where $d_1,\cdots,d_{n-1}$ are domains in $p$}}
\end{array}}
\right\}.
\eeqn
The following proposition guarantees
the existence of $\mpath(n;\vec t;\vec c)$.
\newtheorem{pro7}[def7]{Proposition}
\begin{pro7}
\label{non-emp}
Suppose that $n-|m|\in2\ZZ_{\geq 0}$.
For any $\vec t=(t_1,\cdots,
t_{n-1})$ in $m$-domain configuration and
any sequence of non-negative integers
$\vec c=(c_1,\cdots,c_{n-1})$,
\beq
\mpath(n;\vec t;\vec c)\ne \emptyset
\label{nonempty}
\eeq
\end{pro7}
{\sl Proof.}\,\,
By Lemma \ref{n-m}, if $n-|m|\in2\ZZ_{\geq 0}$,
there exists $\vec t=(t_1,\cdots,t_{n-1})$
in $m$-domain configuration.
Let $p^{(\pm)}_l$ be paths given as follows:
for $j=1,\cdots,n-1$
\beq
d_j:=\left\{
\begin{array}{ll}
\underbrace{\pm t_j,\mp t_j,\cdots,\pm t_j,\mp t_j}_{2c_j}
&{\hbox{ if $t_j$ is  up-regular}},\\
\underbrace{\mp t_j,\pm t_j,\cdots,\mp t_j,\pm t_j}_{2c_j}
&{\hbox{ if $t_j$ is  down-regular}},\\
\underbrace{\pm t_j,\mp t_j,\cdots, \mp t_j,\pm t_j}_{2c_j+1}
&{\hbox{ if $t_j$ is  maximal}}, \\
\underbrace{\mp t_j,\pm t_j,\cdots,\pm t_j,\mp t_j}_{2c_j+1}
&{\hbox{ if $t_j$ is  minimal}},
\end{array}
\right.
\label{tttt}
\eeq

\beq
d_n:=\left\{
\begin{array}{ll}
\pm m,\mp m,\cdots & {\hbox{ if $d_n$ is up,}}\\
\mp m,\pm m,\cdots & {\hbox{ if $d_n$ is down.}}
\end{array}
\right.
\label{mmmm}
\eeq
Now, we order these domains by setting
the position of the left most  $m$ (resp. $-m$)
in $d_n$ being $2l$ (resp. $2l-1$). For example,
$$
p^{(+)}_l=(\cdots 00|d_1|d_2|
\cdots|d_{n-1}|\stackrel{2l}{m}-m\cdots)\,\,{\hbox{ or }}\,\,
(\cdots 00|d_1|d_2|\cdots|d_{n-1}|-m \stackrel{2l}{m}\cdots).
$$
For $p^{(+)}_l$,
by using induction on the index of domains
we shall show the claim
that the position of any entry $t_j$ in $d_j$
is even and the one of $-t_j$ in
$d_j$ is odd. Now we assume that $d_n$ is up.
Then $d_{n-1}$ must be up-regular or minimal
by Definition \ref{def:domaintype}.
It is trivial that in both cases by (\ref{tttt})
the position of $t_{n-1}$ is even and
the one of $-t_{n-1}$ is odd.
Now, we assume that for $i=j+1$ the claim is valid.
If $t_{j+1}$ is up-regular or maximal,
by  Definition \ref{def:domaintype},
$t_j$ must be up-regular or minimal.
Then by (\ref{tttt}) we have
\beq
(\cdots d_j|d_{j+1}\cdots)
=(\cdots t_j,-t_j|t_{j+1},-t_{j+1},\cdots).
\label{733}
\eeq
This implies that the statement is valid for $i=j$.
If $t_{j+1}$ is down-regular or minimal,
by Definition \ref{def:domaintype},
 $t_j$ must be down-regular or maximal.
Then by (\ref{tttt}) we have
\beq
(\cdots d_j|d_{j+1}\cdots)=(\cdots -t_j,t_j|-t_{j+1},t_{j+1},\cdots).
\label{734}
\eeq
This implies that the statement is valid for $i=j$.
Therefore, we have
$$
t_j=t(d_j)\,\,{\hbox{ and then }}\,\, c_j=c(d_j).
$$
We obtain that $p^{(+)}_l\in\mpathtc$.
We can also show for $p^{(-)}_l$. \qed

\vskip5pt
\nd
{\sl Reamrk.}\,
These $p^{(\pm)}_l$ are  extremal vectors in $\mpathtc$,
which is shown in \ref{Ext in mpath}.
We also know that all walls in $p^{(+)}_l$
(resp. $p^{(-)}_l$) are $+$ (resp. $-$)
simultaneously by the definition of $p^{(\pm)}_l$ and
the conditions in Definition \ref{def:domaintype}.
In fact, in (\ref{733}) (resp. (\ref{734})), we get
$t_{j+1}-t_j=1$ and  (resp.  $t_{j}-t_{j+1}=1$).

\subsection{Stability of $\mpath(n;\vec t,\vec c)$}
We shall show the stability of $\mpath(n;\vec t;\vec c)$
by the actions of $\eit$ and $\fit$.

\begin{pro7}
\label{stab:pntc}
For any $i\in I$, we have
\beq
\begin{array}{rcl}
\eit \mpath(n;\vec t;\vec c) & \subset
 & \mpath(n;\vec t;\vec c)\sqcup \{0\},\\
\fit \mpath(n;\vec t;\vec c) & \subset
 & \mpath(n;\vec t;\vec c)\sqcup \{0\}.
\end{array}
\eeq
\end{pro7}

In order to show this proposition,
we shall prepare several lemmas.
\begin{lem7}
\label{even-odd}
For $p=(\cdots,i_{k-1},\,i_k,\,i_{k+1},\cdots)\in \mpath(n)$,
suppose that
$\til f_ip=(\cdots,i_{k-1},\,\til f_i(i_k),\,i_{k+1},\cdots)$
$($resp. $\til e_ip=(\cdots,i_{k-1},
\,\til e_i(i_k),\,i_{k+1},\cdots))$ and let
$d_j$ be the domain including $i_k$. Then we have
\begin{enumerate}
\item
The entry $i_k$ is the right-most entry
$({\hbox{resp. left-most entry}})$
in $d_j$.
\item
Suppose that $d_j$ is a finite domain.
The length $l(d_j)$ is odd if and only if
$d_j$ is regular and the length $l(d_j)$ is even
if and only if $d_j$ is critical.
\item
Suppose that $d_{j+1}$ is a finite domain.
The length $l(d_{j+1})$ $($resp. $l(d_{j-1}))$ is even
if and only if $d_{j+1}$ $($resp. $d_{j-1})$ is regular
and the length $l(d_{j+1})$ $($resp. $l(d_{j-1}))$ is odd
if and only if $d_{j+1}$ $($resp. $d_{j-1}$$)$ is critical.
\end{enumerate}
\end{lem7}
{\sl Remark.}\,
The statement (i) means that $d_j\ne d_n$ (resp. $d_j\ne d_0$), that is ,
there is a domain on the right (resp. left) side of $d_j$.
Then, the statement (iii) makes sense.

\vskip6pt
\nd
{\sl Proof.}\,
Since the proof for the $\til e_i$ case is similar
to the one for $\til f_i$,
we shall show only for the $\til f_i$ case.

\vskip5pt
\nd
(i)\,\,
By Lemma \ref{action} (i),
the hypothesis
$\til f_ip=(\cdots,i_{k-1},\til f_i(i_k),i_{k+1},\cdots)$
implies that $a^{(i)}_k<a^{(i)}_{k+1}$ and then we have
\beqn
&& i_k+i_{k+1}<0\,\, {\hbox{ if $i=1$ }},
\label{kkkk}\\
&& i_k+i_{k+1}>0\,\,{\hbox{ if $i=0$}}.
\eeqn
Then we get
$|i_k+i_{k+1}|>0$. This gives the desired result.

\vskip5pt
\nd
(ii)\,\,
We shall show the $\til f_1$-case.
Let $i_r$ be the left-most entry in $d_j$.
(by (i) the right-most entry is $i_k$, then $r\leq k$.).
We set $t:=t(d_j)$, then, $i_r=\pm t$ and $i_k=\pm t$.
Let us recall $a^{(i)}_k$ in (\ref{a_k}).
Owing to (\ref{l-l}) and $\vp_i(x)-\vep_i(-x)=0$,
we have
$a^{(1)}_r=a^{(1)}_k.$
Then by Lemma \ref{action}, we get
$a^{(1)}_r=a^{(1)}_k\leq a^{(1)}_{r-1}$ and then
\beq
0\leq a^{(1)}_{r-1}-a^{(1)}_r
=-\vp_1(i_{r-1})+\vep_1(i_r)=i_{r-1}+i_r.
\label{ip-1+ip>=0}
\eeq
The definition of $i_r$ that $i_r$ is the left-most entry in $d_j$
implies that there are walls
in position $r$ and then $i_{r-1}+i_r\ne 0$.
Thus, due to (\ref{ip-1+ip>=0}) we get
\beq
i_{r-1}+i_r>0.
\label{ip-1+ip>0}
\eeq
There are the following cases (a)--(d):
\begin{enumerate}
\renewcommand{\labelenumi}{(\alph{enumi})}
\item
$i_r=i_k=t$. ({\it i.e.} $r$ and $k$ are even.)
\item
$i_r=i_k=-t$. ({\it i.e.} $r$ and $k$ are odd.)
\item
$i_r=-t$ and $i_k=t$. ({\it i.e.} $r$ is odd and $k$ is even.)
\item
$i_r=t$ and $i_k=-t$. ({\it i.e.} $r$ is even and $k$ is odd.)
\end{enumerate}

\nd
In fact,
the condition (a) or (b) is equivalent
to that $l(d_j)=k-r+1$ is odd and
the condition (c) or (d) is equivalent
to that $l(d_j)=k-r+1$ is even.
Since these (a)--(d) covers all possibilities for $d_j$,
it is enough to show that
if (a) or (b), $d_j$ is regular and if (c) or (d),
$d_j$ is critical.
Let $d_s$ and $d_p$ be the domains including
$i_{r-1}$ and $i_{k+1}$ respectively.

\vspace{4pt}
\nd
In the case (a) (resp. (b)),
by (\ref{kkkk}) and (\ref{ip-1+ip>0}),
we get $i_{k+1}<-t$ (resp. $i_{k+1}<t$)
and  $i_{r-1}>-t$ (resp. $i_{r-1}>t$).
Since $k+1$ and $r-1$ are odd (resp. even), the domain types
$t(d_p)=- i_{k+1}>t$ (resp. $t(d_p)=i_{k+1}<t$)
and  $t(d_s)=-i_{r-1}<t$ (resp. $t(d_s)=i_{r-1}>t$).
This implies
\beq
t(d_{j+1})=t+1,\,\, t(d_{j-1})=t-1\q
{\hbox{$($resp. $t(d_{j+1})=t-1,\,\, t(d_{j-1})=t+1$.$)$}}.
\label{tdi+1=t+1,tdp-1=t-1}
\eeq
Furthermore, this (\ref{tdi+1=t+1,tdp-1=t-1})
implies that the domain $d_j$ is regular.

\vspace{4pt}
\nd
In the case (c) (resp. (d)),
by (\ref{kkkk}) and (\ref{ip-1+ip>0}),
we get $i_{k+1}<-t$ (resp. $i_{k+1}<t)$
and $i_{r-1}>t$ (resp. $i_{r-1}>-t$).
Since $k+1$ is odd (resp. even) and $r-1$
is even (resp. odd), the domain types
$t(d_p)=-i_{k+1}>t$ (resp. $t(d_p)=i_{k+1}<t$)
and  $t(d_s)=i_{r-1}>t$ (resp.  $t(d_s)=-i_{r-1}<t$).
This implies that
\beq
t(d_{j+1})=t+1,\,\, t(d_{j-1})=t+1\q
{\hbox{$($resp. $t(d_{j+1})=t-1,\,\, t(d_{j-1})=t-1$.$)$}}.
 \label{tdi+1=t+1,tdp-1=t+1}
\eeq
Furthermore, this (\ref{tdi+1=t+1,tdp-1=t+1})
means that the domain $d_j$ is critical.
Now, we have completed the proof of (ii)

\vskip5pt
\nd
(iii)\,\,
We shall show the $\til f_1$-case.
Since  $i_k+i_{k+1}<0$ by (\ref{kkkk}),
we shall consider the following two cases:
\begin{enumerate}
\renewcommand{\labelenumi}{(\arabic{enumi})}
\item
 $i_k+i_{k+1}\leq -2$
\item
 $i_k+i_{k+1}=-1$.
\end{enumerate}

\vskip5pt
\nd
(1)
The assumption $i_k+i_{k+1}\leq -2$ implies that
the domain $d_{j+1}$ is a domain with zero-length.
By Remark under Definition \ref{def:domaintype},
$d_{j+1}$ is a regular domain.

\vskip5pt
\nd
(2)\,\,
The assumption $i_k+i_{k+1}=-1$ means
that $i_{k+1}=\pm t-1$ $(t=t(d_j))$ and
there is only one wall in position $k+1$.
Then we know that
$i_{k+1}$ is included in the domain
$d_{j+1}$ and $i_{k+1}$ is the
left-most entry of $d_{j+1}$.
Let $i_l$ be the right-most entry of $d_{j+1}$
$(k+1\leq l)$.

By the definition of $a^{(i)}_k$, we have
\beq
a^{(1)}_{k+1}  =  a^{(1)}_k+\vp_1(i_k)-\vep_1(i_{k+1})
               =  a^{(1)}_k-(i_k+i_{k+1})=a^{(1)}_k+1.
\label{ak+1=ak+1}
\eeq
By Lemma \ref{action}
if $\nu>k$, then $a^{(1)}_{\nu}>a^{(1)}_k$.
Then by (\ref{ak+1=ak+1})
we have
\beq
a^{(1)}_{\nu}\geq a^{(1)}_{k+1}\,\,\,(\nu\geq k+1).
\label{anu>ak+1}
\eeq
Owing to (\ref{l-l}), we can easily get
\beq
a^{(1)}_{l+1}=a^{(1)}_{k+1}-(i_l+i_{l+1}).
\label{as+1=ak+1-is}
\eeq
The formula (\ref{anu>ak+1}) and (\ref{as+1=ak+1-is}) show,
\beq
i_l+i_{l+1}\leq 0.
\label{is+is+1<=0}
\eeq
Since $i_l$ is the right-most entry in $d_{j+1}$,
there exist walls in position $l+1$.
Therefore, by (\ref{is+is+1<=0}),
\beq
i_l+i_{l+1}<0.
\label{is+is+1<0}
\eeq
As in (ii),
there are the following four cases (a)--(d)
since $i_{k+1}=\pm t-1$:
\begin{enumerate}
\renewcommand{\labelenumi}{(\alph{enumi})}
\item
$i_{k+1}=i_l=t-1$.
({\it i.e.} $i_k=-t$, $k+1$ and $l$ are even.)
\item
$i_{k+1}=i_l=-t-1$.
({\it i.e.} $i_k=t$, $k+1$ and $l$ are odd.)
\item
$i_{k+1}=t-1$ and $i_l=-t+1$.
({\it i.e.} $i_k=-t$, $k+1$ is even and $l$ is odd.)
\item
$i_{k+1}=-t-1$ and $i_l=t+1$.
({\it i.e.} $i_k=t$, $k+1$ is odd and $l$ is even)
\end{enumerate}

\vskip5pt
\nd
The condition (a) or (b) is equivalent
to that $l(d_{j+1})$ is odd and
the condition (c) or (d) is equivalent
to that $l(d_{j+1})$ is even.
Thus, it is enough to show that
if (a) or (b), $d_{j+1}$ is  critical  and if (c) or (d),
$d_{j+1}$  is regular.

\nd
Let $d_q$ be the domain including $i_{l+1}$.
Applying (\ref{is+is+1<0}) to these cases,
we get
\begin{enumerate}
\renewcommand{\labelenumi}{(\alph{enumi})}
\item
$t(d_q)=-i_{l+1}>i_l=t-1=t(d_{l+1})$. This implies that
  $t(d_{j+2})=t$ and then $d_{j+1}$ is a critical domain.
\item
$t(d_q)=i_{l+1}<-i_l=t+1=t(d_{l+1})$. This implies that
  $t(d_{j+2})=t$ and then $d_{j+1}$ is a critical domain.
\item
$t(d_q)=i_{l+1}<-i_l=t-1=t(d_{l+1})$. This implies that
  $t(d_{j+2})=t-2$ and then $d_{j+1}$ is a regular domain.
\item
$t(d_q)=-i_{l+1}>i_l=t+1=t(d_{l+1})$.  This implies that
  $t(d_{j+2})=t+2$ and then $d_{j+1}$ is a regular domain.
\end{enumerate}

\nd
Since the cases (a)--(d) cover
all possibilities for $d_{j+1}$,
we obtain the desired results.\qed

Now, we set that
for domains $d=i_k,\cdots i_l$ $($$k\leq l$$)$
in a path $p$  and $d'=j_s,\cdots,j_t$
$($$s\leq t$$)$ in a path $p'$,
$d\subset d'$ if $s\leq k\leq l\leq t$
and $i_r=j_r$ for $r=k,\cdots,l$. We set
$d=d'$ if and only if $d\subset d'$ and $d'\subset d$.

\begin{lem7}
\label{tc}
Suppose that for
$p=(\cdots,i_{k-1},i_k,i_{k+1},\cdots)\in \mpath(n)$
$\til f_ip=(\cdots,i_{k-1},\til f_i(i_k),i_{k+1},\cdots)$
$($resp. $\til e_ip=(\cdots,i_{k-1},\til e_i(i_k),i_{k+1},\cdots))$ and
let $d_1,\cdots,$$d_{n-1}$ and $d'_1,\cdots,d'_{n-1}$
be the finite domains in
$p$ and $\til f_ip$ $($resp. $\til e_ip$$)$ respectively.
In particular,
let $d_j$ be the domain including $i_k$.
Then, we get
\begin{enumerate}
\item
If $i\ne j,j+1$ $($resp. $i\ne j-1,j$$)$,
 \,\,{\hbox{ then }}\,\,$d_i=d'_i$.
\item
If the domain $d_j$ is finite, as a set $d'_j\subset d_j$ and
$d_j\setminus d'_j=\{i_k\}$  and then
$$
l(d'_j)=l(d_j)-1
$$
\item
If the domain $d_{j+1}$ $($resp. $d_{j-1}$$)$ is finite,
as a set $d_{j+1}\subset d'_{j+1}$
$($resp. $d_{j-1}\subset d'_{j-1}$$)$ and
$d'_{j+1}\setminus d_{j+1}=\{\fit (i_k)\}$
$($resp. $d'_{j-1}\setminus d_{j-1}=\{\eit (i_k)\}$$)$ and then
$$
l(d'_{j+1})=l(d_{j+1})+1
\q({\hbox{resp. }}l(d'_{j-1})=l(d_{j-1})+1).
$$
\end{enumerate}
\end{lem7}
{\sl Proof.}\,
We shall see only the $\til f_1$ case
since other cases can be shown similarly.

\nd
By the proof of Lemma \ref{even-odd}, we know that
$i_k+i_{k+1}<0.$
We can also get
$i_{k-1}+i_k\geq 0$.
By the fact that $\til f_1(i_k)=i_k+1$, we have
\beq
|\til f_1(i_k)+i_{k+1}|=|i_k+i_{k+1}|-1\,\,{\hbox { and }}\,\,
|i_{k-1}+\til f_1(i_k)|=|i_{k-1}+i_k|+1.
\label{f1ik}
\eeq
This means that one wall in position $k+1$ shifts to position $k$ and
the entry in the position $k$ is transferred
from $d_j$ to $d'_{j+1}$ by the action of $\til f_1$.
The shifted  wall is the $j+1$ th wall
 since it is on the right boundary
of the domain $d_j$. Here note that
a domain $d_k$ is surrounded by $k$-th wall
and $k+1$-th wall.
Thus we obtain the desired results.\qed

\vskip5pt
\nd
{\sl Proof of Proposition \ref{stab:pntc}.}\,

\nd
For $p\in\mpathtc$, suppose that
$\til f_ip=(\cdots,i_{k-1},\til f_i(i_k),i_{k+1},\cdots)\ne 0.$
Let $d_0,d_1, \cdots,d_{n-1},d_n$ and $d'_0,d'_1, \cdots,d'_{n-1},d'_n$
be domains
of $p$ and $\til f_1p$ respectively,
in particular $d_j$ be the domain including $i_k$
($d_0$, $d_n$, $d'_0$ and $d'_n$ are infinite domains.).
First let us show
\beq
t(d'_i)=t(d_i)\,\,{\hbox { for any }}\,i=1,2,\cdots,n-1.
\label{td=td}
\eeq
By Lemma \ref{tc} (i),
we know that for $i\ne j,j+1$ such that $d_i=d'_i$ is
non-zero length domain,
$t(d'_i)=t(d_i).$
We shall consider the type of $d'_{j+1}$.
If $d_{j+1}$  and $d'_{j+1}$ are  infinite domains,
$j+1=0$ or $n$ then there is nothing to prove.
Then we may assume that $d_{j+1}$ and $d'_{j+1}$ are finite domains.
If $l(d_{j+1})\geq 1$,
there exists $a\in\ZZ$ such that $i_a$ is
included in both $d_{j+1}$ and $d'_{j+1}$ by Lemma \ref{tc} (iii).
  Then, in this case we get
$t(d'_{j+1})=t(d_{j+1})$.
In the case $l(d_{j+1})=0$
if we assume that $t(d_j)=t$ and $t(d_{j+1})=t+1$,
by  the proof of Lemma \ref{even-odd} related to (a) and (c),
we get that $i_k=t$ and $k$ is even.
Then $\til f_1(i_k)=t+1$.
This entry is included in $d'_{j+1}$ and then
$t(d'_{j+1})=(-1)^k\til f_1(i_k)=(-1)^k(t+1)=t+1$.
We can also easily see the case $t(d_{j+1})=t-1$.
Thus we get
$$
t(d'_{j+1})=t(d_{j+1}).
$$
We shall consider the type of $d'_j$.
As same as above,
 we may assume that $d_j$ and $d'_j$ are finite domains.
If $l(d_j)\geq 2$,
there exists $b\in\ZZ$ such that
$i_b$ is included in both $d_{j}$ and $d'_{j}$
 by Lemma \ref{tc} (ii).
  Then, in this case we get
$t(d'_{j})=t(d_{j})$.
If $l(d_j)=1$, by Lemma \ref{even-odd} (ii)
and Lemma \ref{tc} (ii),  we get that
$d_j$ is a regular domain and
$l(d'_j)=0.$
Since by the previous arguments
we have already obtained $t(d'_i)=t(d_i)$ for $i\ne j$
such that $d_i$ or $d'_i$ is a non-zero length domain and
$d'_j$ is a regular domain
by the remark under Definition \ref{def:domaintype},
 we get
$$
t(d'_j)=t(d_j).
$$
Thus we get $t(d_i)=t(d'_i)$ for all other zero-length domains.
Then we obtain (\ref{td=td}).

Next,  let us show
\beq
c(d'_i)=c(d_i)\,\,{\hbox{ for any }}\,i=1,2,\cdots,n-1.
\label{cd=cd}
\eeq
By (\ref{td=td}), $d_i$ is a regular (resp. critical) domain
if and only if $d'_i$ is a regular (resp. critical) domain.
Therefore, by Lemma \ref{tc} (i)
we have
\beq
c(d'_i)=c(d_i)\,\,{\hbox{ for }}\,\, i\ne j,j+1.
\label{cd'=cd}
\eeq
We shall consider the domain parameter $c(d'_{j})$.
We may assume that $d_{j}$ and $d'_{j}$ are finite domains as in the
previous arguments.
If $d_j$  is a regular, $d'_j$ is also regular
and by Lemma \ref{even-odd} (ii)
and Lemma  \ref{tc} (ii)
we have
\beq
l(d_j)  =  2c_j+1 \,\,{\hbox{ and }}\,\,
l(d'_j)  =  2c_j.
\label{ld'j=2c}
\eeq
Since $d_j$ and $d'_j$ are  regular domains,
the formula (\ref{ld'j=2c}) implies
\beq
c(d'_j)=c_j=c(d_j).
\label{cd'j=c}
\eeq
If $d_j$  is a critical domain, $d'_j$ is also critical
 and by Lemma \ref{even-odd} (ii)
and Lemma  \ref{tc} (ii)
we have
\beq
l(d_j)  =  2c_j+2
 \,\,{\hbox{ and }}\,\,
l(d'_j)  =  2c_j+1.
\label{ld'j=2c+1}
\eeq
Since $d_j$ and $d'_j$ are critical domains,
the formula (\ref{ld'j=2c+1}) implies
\beq
c(d'_j)=c_j=c(d_j).
\label{cd'j=c:2}
\eeq
As for $d'_{j+1}$, by using
Lemma
\ref{even-odd} (iii)
and
Lemma \ref{tc} (iii) we can also easily obtain
\beq
c(d'_{j+1})=c_{j+1}=c(d_{j+1}).
\label{djjj}
\eeq
Thus by (\ref{cd'=cd}), (\ref{cd'j=c})  (\ref{cd'j=c:2}) and
(\ref{djjj}) we get (\ref{cd=cd}).
Now, we have completed
the proof of Proposition \ref{stab:pntc}.\qed

\subsection{Extremal vectors in $\mpath(n;\vec t;\vec c)$}
\label{Ext in mpath}
In this subsection,
we shall describe all extremal vectors in $\mpath(n;\vec t;\vec c)$
explicitly.

\begin{lem7}
\label{lem:ext}
Let $ B_1$ and $B_2$ be normal crystals
and $\phi:B_1\rightarrow B_2$ be a
strict morphism of crystal and we assume
that $\phi(b)\ne 0$ for $b\ne 0$.
We have that
$b$ is an extremal vector in $B_1$ if and only if
$\phi(b)$ is an extremal vector in $B_2$.
\end{lem7}
{\sl Proof.\,}
We assume that $b$ is not an extremal vector in $B_1$ and
$\phi(b)$ is an extremal vector in $B_2$.
Then there exist $i,i_1,\cdots,i_k\in I$ such that
$$
\eit S_{i_1}\cdots S_{i_k}b \ne 0\,\,{\hbox{ and }}\,\,
\fit S_{i_1}\cdots S_{i_k}b \ne 0.
$$
By the assumption that $\phi(b)\ne 0$ for $b\ne 0$,
we get $\phi(\eit S_{i_1}\cdots S_{i_k}b) \ne 0$ and
$\phi(\fit S_{i_1}\cdots S_{i_k}b) \ne 0$.
Since $\phi$ is a morphism of crystal, we have
$$
\eit S_{i_1}\cdots S_{i_k}\phi(b) \ne 0\,\,{\hbox{ and }}\,\,
\fit S_{i_1}\cdots S_{i_k}\phi(b) \ne 0.
$$
This contradicts  the fact that $\phi(b)$ is an extremal vector.
If $b$ is an extremal vector in $B_1$
and $\phi(b)$ is not an extremal vector in $B_2$,
by arguing similarly we can obtain a contradiction. \qed

Let $\vec t=(t_1,\cdots,t_{n-1})$ be
a $m$-domain configuration($m\in\ZZ$) and
$\vec c=(c_1,\cdots,c_{n-1})$ be
a sequence of non-negative integers.
For $p\in \mpathtc\ne \emptyset$,
let $d_1,\cdots, d_{n-1}$ be its
finite domains. For a domain $d_i$ we set
\beq
l_{\rm min}(d_j):=\left\{
\begin{array}{ll}
2c_i    &{\hbox{ if $d_i$ is regular}}\\
2c_i+1  &{\hbox{ if $d_i$ is critical}}
\end{array}
\right.
\label{lmin}
\eeq
\begin{lem7}
\label{min+-}
For $p\in \mpath(n)$ let
$\io_1(p),\cdots,\io_n(p)$ be its types of walls and set
\beqn
A_m(n;\vec t;\vec c) & := &
\left\{p\in\mpathtc|
\io_1(p)=\cdots=\io_n(p)\right\}.
\label{AM}\\
E_m(n;\vec t;\vec c) & := &
\left\{p\in\mpathtc|l(d_i)=l_{\rm min}(d_i)
\,\,{\hbox{ for any $i$.}}\right\}.
\label{EM}
\eeqn
Then we get
\beq
A_m(n;\vec t;\vec c)=E_m(n;\vec t;\vec c).
\label{A=M}
\eeq
\end{lem7}
{\sl Proof.}\,
For $p\in E_m(n;\vec t;\vec c)$ suppose that
a domain $d_j$  in $p$ is a regular domain with non-zero length
and set $t(d_j)=t$.
Let $i_a$ and $i_b$ be the left-most entry
and the right-most entry in $d_j$ respectively. By (\ref{lmin}),
$l(d_j)=b-a+1=2c_j>0$.
Thus, if $a$ is even (resp. odd), $b$ is odd (resp. even).
Now we assume that $a$ is even and $b$ is odd. Then we have
\beq
t(d_j)=i_a=-i_b.
\label{tdj=ia}
\eeq
Let $d_r$ and $d_s$ be
the domains including $i_{a-1}$ and $i_{b+1}$ respectively.
We have
\beq
t(d_r)=-i_{a-1}\,\,{\hbox{ and }}\,\,t(d_s)=i_{b+1},
\label{tdr,tds}
\eeq
since $a-1$ is odd and $b+1$ is even.
Because $d_j$ is regular,
\beq
t(d_r)<t(d_j)<t(d_s)\,\,{\hbox{ or }}\,\,
t(d_r)>t(d_j)>t(d_s).
\label{or}
\eeq
Applying (\ref{tdj=ia}) and
(\ref{tdr,tds}) to (\ref{or}) we obtain
$$
i_{a-1}+i_a>0,\,i_b+i_{b+1}>0\,\,{\hbox{ or }}\,\,
i_{a-1}+i_a<0,\,i_b+i_{b+1}<0.
$$
This means that all walls in $a$
and in $b+1$  have the same type.
We can get the same result for the case
that $a$ is odd and $b$ is even,
and the case that $d_j$ is critical.
Repeating this for all domains with non-zero length,
we know that all walls have same type in $p$. Thus, we  have
\beq
E_m(n;\vec t;\vec c)\subset A_m(n;\vec t;\vec c).
\label{E<A}
\eeq
Let $p$ be an element of $A_m(n;\vec t;\vec c)$
and all walls in $p$ be $+$.
For a  regular domain with non-zero length $d_j$ in $p$
let $i_a$ and $i_b$ be
left-most entry and right-most entry in $d_j$ respectively.
Then we get
\beq
i_{a-1}+i_a>0,\,{\hbox{ and }}\,i_b+i_{b+1}>0.
\label{aabb}
\eeq
Let $d_r$ and $d_s$ be as above.
Since $d_j$ is regular, we have (\ref{or}).
If $a$ is even, $t(d_j)=i_a$ and $t(d_r)=-i_{a-1}$.
By (\ref{aabb}), we get $t(d_r)<t(d_j)$.
Thus, by the assumption  that $d_j$ is regular, we have
\beq
t(d_r)<t(d_j)<t(d_s).
\label{rjs}
\eeq
Furthermore, if $b$ is even, $t(d_j)=i_b$ and $t(d_s)=-i_{b+1}$.
Then this and
(\ref{rjs}) imply that $i_b+i_{b+1}<0$.
But this contradicts  (\ref{aabb}).
Then $b$ is odd and then $l(d_j)=b-a+1$ is even.
Since $d_j$ is regular, this means
$$
l(d_j)=2c_j=l_{\rm min}(d_j).
$$
By arguing similarly for other non-zero length domains, we obtain
$l(d_i)=l_{\rm min}(d_i)$ for any $i$. Therefore, we get
\beq
A_m(n;\vec t;\vec c) \subset E_m(n;\vec t;\vec c).
\label{A<E}
\eeq
By (\ref{E<A}) and (\ref{A<E}), we get the desired result.\qed

\begin{pro7}
\label{extvec}
Let $E$ be the set of all extremal vectors in $\mpathtc$.
Then we have
\beq
E=A_m(n;\vec t;\vec c)=E_m(n;\vec t;\vec c).
\label{E=E}
\eeq
\end{pro7}
{\sl Proof.}\,
By the definition of the map $\psi$ given in (\ref{V}),
we know that $\psi(p)\ne 0$ for $p\in\mpath(n)$. Therefore,
by Proposition \ref{ext}, Theorem \ref{path-spin}
and Lemma \ref{lem:ext},
we get
\beq
 A_m(n;\vec t;\vec c)= E.
\eeq
By Lemma \ref{min+-}, we obtain the desired result. \qed

Let $p^{(\pm)}_l$ be paths given
in the proof of Proposition \ref{non-emp} and
set
$$
E':=\{p^{(\pm)}_l\}_{l\in \ZZ}.
$$
\newtheorem{thm7}[def7]{Theorem}
\begin{thm7}
We have
\beq
E=E'=A_m(n;\vec t;\vec c)=E_m(n;\vec t;\vec c).
\label{E=E'}
\eeq
\end{thm7}
{\sl Proof.}\,\,
By the definiton of $p^{(\pm)}_l$
in the proof of Proposition \ref{non-emp} and
Proposition \ref{extvec}, we get $E'\subset E$ easily.
For $p^{(\ep)}_l$ ($\ep=\pm$, $l\in\ZZ$)
let $k^{\ep,l}_1,\cdots,k^{\ep,l}_n$ be the
positions of walls in $p^{(\ep)}_l$.
By (\ref{tttt}), (\ref{mmmm}) and the way of ordering,
we get
\beq
\begin{array}{rcl}
(k^{+,l}_1,\cdots,k^{+,l}_n)
& = & (k^{-,l}_1+1,\cdots,k^{-,l}_n+1),\\
(k^{-,l}_1,\cdots,k^{-,l}_n)
& = & (k^{+,l-1}_1+1,\cdots,k^{+,l-1}_n+1).
\end{array}
\label{kkkk-2}
\eeq
Let $p$ be an element in $E_m(n;\vec t;\vec c)$ and
$(k_1,\cdots,k_n)$ be the positions of walls in $p$.
By the definiton of $E_m(n;\vec t;\vec c)$,
we know that for any $\ep$ and $l$,
$ k^{\ep,l}_{j+1}-k^{\ep,l}_j=l_{\rm min}(d_j)=(2c_j$
or $2c_j+1)=k_{j+1}-k_j$.
Therefore, by (\ref{kkkk-2}),
there exist $\ep\in\{\pm\}$
and $l\in\ZZ$ such that
$(k_1,\cdots,k_n)=(k^{\ep,l}_1,\cdots,k^{\ep,l}_n)$.
Now, since the domain types are fixed,
the entries in $p$ are automatically determined and
it coincides with the ones in $p^{(\ep)}_l$.
This means that $p=p^{(\ep)}_l$ and then
$$
E_m(n;\vec t;\vec c)\subset E'.
$$
Now, we have completed the proof.\qed

\vskip7pt
\nd
{\sl Remark.}\,
By Lemma \ref{tc}, we know that
a $(-)$(resp. $(+)$) wall in a path is shifted
by one to the left direction by the action of
$\til f_1$ (resp. $\til f_0$) and
by the definition of $p^{(\pm)}_l$
$(k^+_1,\cdots,k^+_n)=(k^-_1+1,\cdots,k^-_n+1)$,
where $(k^{\pm}_1,\cdots,k^{\pm}_n)$ are
sequences of the positions of walls in $p^{(\pm)}_l$.
Therefore, we have
\beq
\til f_1^np^{(-)}_l=p^{(+)}_{l-1}\,\,{\hbox{ and }}\,\,
\til f_0^np^{(+)}_l=p^{(-)}_{l}.
\label{f1f0}
\eeq
Thus, we have
\beq
S_1p^{(-)}_l=\til f_1^np^{(-)}_l=p^{(+)}_{l-1}\,\,
{\hbox{ and }}\,\,S_0p^{(+)}_l=\til f_0^np^{(+)}_l=p^{(-)}_l.
\label{S1fS0f}
\eeq
By these (\ref{f1f0}) and (\ref{S1fS0f}), we get
\beq
S_1p^{(-)}_l=p^{(+)}_{l-1,}\,\,
S_0p^{(+)}_l=p^{(-)}_{l},\,\,
S_1p^{(+)}_l=p^{(-)}_{l+1}
{\hbox{ and }}\,\,
S_0p^{(-)}_l=p^{(+)}_l.
\label{weyl}
\eeq
Thus, we obtain the following result.
\newtheorem{cor7}[def7]{Corollary}
\begin{cor7}
\label{connected}
$\mpathtc$ is a connected component in $\mpath$.
\end{cor7}
{\sl Proof.}\,
By the remark as above, we know that $E'$ is connected and then any extremal
vector
in $\mpathtc$ is connected to each other.
Therefore, by Theorem \ref{conn-ext} and Proposition \ref{stab:pntc},
we know that $\mpathtc$ is connected.\qed

\begin{ex7}
$$
\begin{array}{llllllllllll}
&\cdots &0&0&0\stackrel{+}{|}
&1&-1\stackrel{+}{|}&2&-2\stackrel{+}{|}&3&-3&\cdots\\
\mapleftright{S_0}
&\cdots &0&0\stackrel{-}{|}
&-1&1\stackrel{-}{|}&-2&2\stackrel{-}{|}&-3&3&-3&\cdots \\
\mapleftright{S_1}
&\cdots &0\stackrel{+}{|}
&1&-1\stackrel{+}{|}&2&-2\stackrel{+}{|}&3&-3&3&-3&\cdots
\end{array}
$$
\end{ex7}

\subsection{Affinization of the Path-Spin Correspondence}
In Sec.6 we introduced the path-spin correspondence. In this subsection,
we shall affinize it, that is, the path-spin correspondence in Sec.6,
which  is a morhpism of classical crystal,
is lifted to a morphism of affine crystals.

\vskip5pt
\nd
Let $B=\{+,-\}$ be the classical crystal as in Example \ref{ex4:B}.
\begin{lem7}
\label{lem:zn} $(${\rm see 2.2}$)$\,\,
The set of all extremal vectors in $\aff(B^{\ot n})$ is given by
\beq
\{z^k\ot (+)^{\ot n},z^k\ot (-)^{\ot n}\}_{k\in \ZZ}.
\nonumber
\eeq
\end{lem7}
{\sl Proof.}\,
By (\ref{delta}), (\ref{act}) and (\ref{eqn:SS}), we have
\beqnn
S_1(z^k\ot (\pm)^{\ot n}) & = & z^k\ot (\mp)^{\ot n},
\\
S_0(z^k\ot (\pm)^{\ot n}) & = & z^{k\pm n}\ot (\mp)^{\ot n}.
\eeqnn
By (\ref{act}) and (\ref{4act}), we get for any $k$
\beq
\til e_1(z^k\ot (+)^{\ot n})=\til f_0(z^k\ot (+)^{\ot n})=
\til e_0(z^k\ot (-)^{\ot n})=\til f_1(z^k\ot (-)^{\ot n})=0.
\nonumber
\eeq
Thus, we get the desired result.\qed

 Now we shall consider the affinization of the morphism $\psi$.
 For a level 0 affine weight
$\lm=m(\Lm_0-\Lm_1)+l\del\in P$ ($l,m\in \ZZ$)
 by Remark (ii) in 4.2,
we have that $Ua_{\lm}$ has a $U'$-module structure
 and its crystal $B(Ua_{\lm})$ is
described by $\mpath$ as a classical crystal
 (that is, $B(Ua_{\lm})\cong B(U'a_{cl(\lm)})$
as a classical crystal).
 Originally,
the crystal $B(Ua_{\lm})$ holds an affine crystal structure.
 We shall recover its affine crystal structure in terms of path.
 For this purpose we shall introduce the {\it energy function}
 (See \cite{IIJMNT}, \cite{KMN}, \cite{KKM}).
\begin{def7}[\cite{KKM}]
 Let B be a classical crystal. A $\ZZ$-valued function $H$ on
 $B\ot B$ is called an {\it energy function} on $B$
if for any $i\in I$ and
 $b\ot b'\in B\ot B$ such that $\eit(b\ot b')\ne 0$, we have
\beqnn
H(\eit(b\ot b'))&=& H(b\ot b'){\hbox{ if }}i=1,\nonumber \\
                &=& H(b\ot b')+1{\hbox{ if }}i=0
                   {\hbox{ and }}\vp_0(b)\geq\vep_0(b'), \\
                &=& H(b\ot b')-1{\hbox{ if }}i=0
                {\hbox{ and }}\vp_0(b)<\vep_0(b').\nonumber
\eeqnn
\end{def7}

\vskip5pt
   For the case of $B=B_{\infty}$, by \cite{KKM} Theorem 5.1,
   we can describe
   the function $H$ explicitly as follows.
\begin{pro7}
\label{Hmax}
We set
$$
  H((m)\ot (n)):={\rm max}\{m,-n\}.
$$
This $H$ is an energy function on $B_{\infty}$.
\end{pro7}

  By Theorem 4.9 in \cite{KKM}, we get the following theorem.

\begin{thm7}
  Let $(g_i)_{i\in \ZZ}$ be a $m$-ground-state path.
For a level 0 affine weight
$\lm=m(\Lm_0-\Lm_1)+l\del\in P$ and $b\in B(Ua_{\lm})$
which corresponds to the $m$-path $p=(i_k)_{k\in \ZZ}\in \mpath$
as a classical crystal, we have the following formula,
\beq
\begin{array}{rcl}
\lefteqn{ wt(b)=wt(p)=(\sum_{k\in\ZZ}i_{k-1}+i_{k})(\Lm_0-\Lm_1) } \\
&& \qq  +(l+\sum_{k\in\ZZ}k({\rm max}\{i_{k-1},-i_{k}\}
-{\rm max}\{g_{k-1}, -g_k\}))\del.
 \end{array}
  \label{path-wt}
 \eeq
\end{thm7}
{\sl Proof.}\,
   By using the formula in \cite[Theorem 4.9.]{KKM}
and the same type of the formula for
   $B(-\infty)$, we can easily derive
\beq
\begin{array}{rcl}
  \lefteqn{ wt(b)=wt(p)=m(\Lm_0-\Lm_1)
+\sum_{k\in\ZZ}(af(wt(i_k))-af(wt(g_k))) } \\
   && \qq+\{l+\sum_{k\in\ZZ}k(H(i_{k-1}
\ot i_{k})-H(g_{k-1}\ot g_k))\}\del.
 \end{array}
  \label{path-wt2}
 \eeq
Since by the definition of path,
the summations in (\ref{path-wt2}) are finite, we get
\beq
 \begin{array}{rcl}
 \lefteqn{    \hspace{-30mm}wt(p)=m(\Lm_0-\Lm_1)
+{\frac{1}{2}}\sum_{i\in\ZZ}
(af(wt(i_{k-1}))+af(wt(i_{k})) }\qq\qq \\

-af(wt(g_{k-1}))-af(wt(g_{k})))&&
\hspace{-7mm} +\{l+\sum_{k\in\ZZ}
k(H(i_{k-1}\ot i_{k})-H(g_{k-1}\ot g_k))\}\del.
\end{array}
\label{www}
\eeq
By applying Proposition \ref{Hmax} and
\beqnn
&&af(wt(i_{k-1}))+af(wt(i_{k}))
=2(\Lm_0-\Lm_1)(i_{k-1}+i_{k}),\\
&&af(wt(g_{k-1}))+af(wt(g_{k}))
=\left\{
\begin{array}{ll}
2m(\Lm_0-\Lm_1)& k=0,\\
0&k\ne0,
\end{array}
\right.\\
\eeqnn
to (\ref{www}),  we obtain the desired result.      \qed

For a level 0 weight $\lm=m(\Lm_0-\Lm_1)+l\del$,
we denote $\mlpath$ for
a set of path corresponding to an element of $B(Ua_{\lm})$,
{\it i.e.} as a set $\mlpath$
is equal to $\mpath$ and a weight is given by (\ref{path-wt}).

By using this formula, we get $\widehat\psi$:
the affinization of $\psi$ as follows:
For $p\in \mlpath$ a map $\widehat\psi$ is given by
\beq
\begin{array}{rcl}
\widehat\psi  :  \mlpath
& \longrightarrow &\aff(B^{\ot n})\\
   p &   \mapsto  & z^{\lan d,wt(p)\ran}\ot \psi(p),
\end{array}
\label{hatpsi}
\eeq
Let us denote also $\til\psi$ for
the restriction of $\widehat\psi$ to
$\mlpathtc$, where $\mlpathtc$ is equal
to $\mpathtc$ as a set and
a weight is given by (\ref{path-wt}).
\begin{thm7}
\label{map:psi}
\begin{enumerate}
\item
The map $\widehat\psi$ and $\til \psi$
are strict morphisms of affine crystals.
\item
The map $\til\psi$ is an injective morphism of affine crystal.
\end{enumerate}
\end{thm7}
{\sl Proof.}\,
(i) \,
It is trivial by (\ref{delta}),
Theorem \ref{path-spin} and Proposition \ref{stab:pntc}.

\vskip6pt
\nd
(ii)\,
In order to show (ii) we shall see the following lemmas :
\begin{lem7}
\label{EEE}
Let $E$ be the set of all extremal vectors in $\mlpathtc$.
 If the map $\til\psi|_E$ is injective,
the map $\til\psi$ is injective.
\end{lem7}
{\sl Proof.}\,
We assume that $\til\psi$ is not injective.
Then there exist $p_1,p_2\in\mlpath$ such that
$p_1\ne p_2$ and $\til\psi(p_1)=\til\psi(p_2)$.
We set $b^*:=\til\psi(p_1)=\til\psi(p_2)\in \aff(B^{\ot n})$.
Due to the connectedness of $B^{\ot n}$,
for this $b^*$ there exist $\til x_{i_1},
\cdots,\til x_{i_l}\in\{\eit,\fit\}_{i=0,1}$
and an extremal vector $v\in \aff(B^{\ot n})$ such that
\beq
v=\til x_{i_1}\cdots \til x_{i_l}(b^*).
\label{vXp}
\eeq
Since $v\ne0$ is an extremal vector,
by Theorem \ref{map:psi} (i) and Lemma \ref{lem:ext},
we have that both $\til x_{i_1}\cdots\til x_{i_l}p_1\ne0$ and
$\til x_{i_1}\cdots\til x_{i_l}p_2\ne0$  are elements in $E$.
The injectivity of $\til\psi|_E$ means
$$
\til\psi(\til x_{i_1}\cdots\til x_{i_l}p_1)
\ne\til\psi(\til x_{i_1}\cdots\til x_{i_l}p_2).
$$
since $p_1\ne p_2$ and then
$\til x_{i_1}\cdots\til x_{i_l}p_1
\ne \til x_{i_1}\cdots\til x_{i_l}p_2$.
But this contradicts the fact that
$$
\til x_{i_1}\cdots \til x_{i_l}\til\psi(p_1)
= v=\til x_{i_1}\cdots \til x_{i_l}\til\psi(p_2).
$$
We have completed the proof of Lemma \ref{EEE}. \qed

\vskip 7pt
\nd
{\sl Proof of Theorem \ref{map:psi}} (ii)\,\,
For a path $p\in \mpath$ let $\io_1(p),\cdots,\io_n(p)$
 be a sequence of the types of the walls in $p$.
We set
\beqn
E_+ & := &
 \{p\in E=E_m(n;\vec t;\vec c)|\io_i(p)=+,\,\,i=1,\cdots,n\},
\label{E+}\\
E_- & := &
 \{p\in E=E_m(n;\vec t;\vec c)|\io_i(p)=-,\,\,i=1,\cdots,n\}.
\label{E-}
\eeqn
These $E_{\pm}$ coincides with
$\{p^{(\pm)}_l\}_{l\in \ZZ}$ respectively.
By (\ref{S1fS0f}) and (\ref{weyl}),
we have the following;
\begin{lem7}
\label{lem:ddd}
For any $p^{(\ep_1)}_k\ne p^{(\ep_2)}_l$
$($$\ep_1,\ep_2=\pm$ and $k,l\in\ZZ$$)$
we have
\beq
wt(p^{(\ep_1)}_k)\ne wt(p^{(\ep_2)}_l).
\label{ddd}
\eeq
\end{lem7}
{\sl Proof.}\,\,
If $\ep_1\ne \ep_2$, $wt(p^{(\ep_1)}_k)\ne wt(p^{(\ep_2)}_l)$
since $wt(p^{(-)}_k)=n(\Lm_0-\Lm_1)+D_1\del$ and
$wt(p^{(+)}_l)=n(\Lm_1-\Lm_0)+D_2\del$
where $D_1$ and $D_2$ are some
integers.
Then we may assume that $\ep_1=\ep_2$.
We set $\ep_1=\ep_2=+$ and $k<l$.
By (\ref{weyl}), we have
$S_0S_1p^{(+)}_l=p^{(+)}_{l-1}.$
This means
$$
(S_1S_0)^{l-k}p^{(+)}_l=p^{(+)}_k.
$$
Since $S_1S_0={\til f_1}^n{\til f_0}^n$ for $p^{(+)}_l$,
we get
\beq
\lan d,wt(p^{(+)}_l)\ran-\lan d,wt(p^{(+)}_k)\ran=(l-k)n>0.
\label{l-kn}
\eeq
Now, we have completed the proof of Lemma \ref{lem:ddd}.\qed

This lemma implies that any extremal vector in $E$
has different weight each other.
Since the morphism of affine crystal $\til \psi$ preserves weight,
now we obtain the injectivity of the map $\til\psi|_E$.
Therefore, by Lemma \ref{EEE}, we get the injectivity of $\til\psi$.
We have completed the proof of Theorem \ref{map:psi}.\qed

\vskip7pt
By the formula $S_1p^{(-)}_l=p^{(+)}_{l-1}$
and $S_1p^{(+)}_l=p^{(-)}_{l+1}$ in (\ref{weyl}),
we get
$$
\lan d, p^{(-)}_l\ran=\lan d,p^{(+)}_{l-1}\ran.
$$
By this and (\ref{l-kn}),
for any extremal vectors $p_1,p_2\in\mlpathtc$, we have
$$
\lan d,wt(p_1)\ran\equiv\lan d,wt(p_2)\ran  \pmod{n}.
$$
By this formula, we obtain the following;

\begin{cor7}
\label{final}
\begin{enumerate}
\item
Set $I_n:=\{0,1,\cdots,n-1\}$ and let
$E_{m,l}(n;\vec t;\vec c)$ be
the set of all extremal vectors in $\mlpathtc$.
Then there exists unique $i\in I_n$ such that
$$
\til\psi(E_{m,l}(n;\vec t;\vec c))
=\{z^{i+kn}\ot(\pm)^{\ot n}\}_{k\in\ZZ}.
$$
\item
Let us denote
 $\aff(B^{\ot n})_i$
for a connected component of $\aff(B^{\ot n})$ generated by
extremal vectors $\{z^{i+kn}\ot(\pm)^{\ot n}\}_{k\in\ZZ}$.
Then as a morphism of affine crystals,
$$
\til\psi:\mlpathtc\mapright{\sim} \aff(B^{\ot n})_i.
$$
\end{enumerate}
\end{cor7}

Now, we shall summarize the classification of paths in
$\mlpath\cong B(\uq a_{\lm})$
($\lm=m(\Lm_0-\Lm_1)+l\del$).
By Corollary \ref{final} (ii),
if we fix one connected component in
$\mlpath(n)$,
each element in the component is
classified by $\aff(B^{\ot n})_i$.
Since $\aff(B^{\ot n})_i$ is generated by
$\{z^{i+kn}\ot (\pm)^{\ot n}\}_{k\in\ZZ}$,
any element in $\aff(B^{\ot n})_i$ is in the following form:
\beq
z^{i_{\io_1,\cdots,\io_{n-1},l}+kn}\ot (\io_1)\ot
 \cdots\ot (\io_{n-1}),
\label{ccc}
\eeq
where $k$ is an integer called {\it depth parameter}
 and $i_{\io_1,\cdots,\io_{n-1},l}\in I_n$ is
determined only by $\io_1,\cdots,\io_{n-1},l$
(if $(\io_1,\cdots,\io_{n-1})=(\pm,\cdots,\pm)$,
$i_{\io_1,\cdots,\io_{n-1},l}=i$.).

Therefore,
for given $m,l\in \ZZ$, by the following parameters:
\beqnn
&&{\hbox{$n\in\ZZ_{\geq 0}$ with $n-|m|\in2\ZZ_{\geq 0}$
(the total number of walls)}},\\
&&{\hbox{$(t_1,\cdots,t_{n-1})$ in $m$-domain configuration
(domain types)}},\\
&&{\hbox{$(c_1,\cdots,c_{n-1})\in \ZZ_{\geq 0}^n$
(domain parameters)}},\\
&&{\hbox{$(\io_1,\cdots,\io_{n-1})$ ($\io_j=\pm$)
(types of walls)}},\\
&&{\hbox{$k\in \ZZ$ (depth parameter)}},
\eeqnn
every path in $\mlpath$ is uniquely classified.

 \end{document}